\documentclass[12pt]{amsart}

\usepackage{amssymb}

\usepackage{graphicx}

\usepackage[cmtip,all]{xy}

\newtheorem{theorem}{Theorem}[section]
\newtheorem{lemma}[theorem]{Lemma}

\theoremstyle{definition}

\newtheorem{example}[theorem]{Example}

\theoremstyle{remark}
\newtheorem{remark}[theorem]{Remark}

\numberwithin{equation}{section}

\usepackage[pdftex]{pict2e}
\usepackage{stmaryrd}
\usepackage{amssymb}
\usepackage{amsfonts}
\usepackage{amsthm}
\usepackage{amsmath}
\usepackage{amscd}
\usepackage[latin2]{inputenc}
\usepackage{t1enc}
\usepackage[mathscr]{eucal}
\usepackage{indentfirst}
\usepackage{graphicx}
\usepackage{graphics}
\usepackage{pict2e}
\usepackage{epic}
\numberwithin{equation}{section}
\usepackage[margin=2.5cm]{geometry}
\usepackage{epstopdf} 
\usepackage{enumerate}

\usepackage{tikz-cd}
\usepackage{mathrsfs}
\usepackage{tikz}
\usetikzlibrary{decorations.markings,positioning}
\usepackage{graphicx}

\newcommand{\abs}[1]{\lvert#1\rvert}
\newcommand{\D}{\mathfrak{D}}
\newcommand{\cc}{\mathfrak{c}}
\newcommand{\p}{\mathcal{P}}
\newcommand{\A}{\mathcal{A}}
\newcommand{\X}{\mathcal{X}}
\newcommand{\T}{\mathcal{T}}
\newcommand{\h}{\mathcal{H}}
\newcommand{\K}{\mathbb{K}}
\newcommand{\Q}{\mathbb{Q}}
\newcommand{\Pp}{\mathbb{P}}
\newcommand{\q}{\mathcal{Q}}
\newcommand{\C}{\mathbb{C}}
\newcommand{\Z}{\mathbb{Z}}
\newcommand{\R}{\mathbb{R}}
\newcommand{\Ll}{\mathbb{L}}
\newcommand{\U}{\mathbb{U}}
\newcommand{\N}{\mathbf{N}}
\newcommand{\M}{\mathcal{M}}
\newcommand{\ol}{\overline}
\newcommand{\tx}{\text}
\newcommand{\ra}{\rightarrow}
\newcommand{\fm}{\mathbb{FM}}

\numberwithin{equation}{section}

\begin{document}

\title[Crepant Transformation correspondence For Toric Stack Bundles]{Crepant Transformation correspondence\\
For Toric Stack Bundles}


\author{Qian Chao}
\address{Department of Mathematics\\ Ohio State University\\ 100 Math Tower, 231 West 18th Ave. \\ Columbus,  OH 43210\\ USA}
\curraddr{School of Mathematics and Statistics\\ Wuhan University\\Wuhan 430072\\ China}
\email{chao.191@buckeyemail.osu.edu, chao818@whu.edu.cn}

\author{Jiun-Cheng Chen}
\address{Department of Mathematics\\ Third General Building\\ National Tsing Hua University\\ No. 101 Sec 2 Kuang Fu Road\\ Hsinchu, 30043\\ Taiwan}
\email{jcchen@math.nthu.edu.tw}

\author[Hsian-Hua Tseng]{Hsian-Hua Tseng}
\address{Department of Mathematics\\ Ohio State University\\ 100 Math Tower, 231 West 18th Ave. \\ Columbus,  OH 43210\\ USA}
\email{hhtseng@math.ohio-state.edu}
\thanks{}


\date{\today}

 
\begin{abstract}
We prove a crepant transformation correspondence in genus zero Gromov-Witten theory for toric stack bundles related by crepant wall-crossings of the toric fibers. Specifically, we construct a symplectic transformation that identifies $I$-functions of toric stack bundles suitably analytically continued using Mellin-Barnes type integrals. We compare our symplectic transformation with a Fourier-Mukai isomorphism between the $K$-groups.
\end{abstract}

\maketitle
\tableofcontents
\setcounter{section}{-1}

\section{Introduction}
\subsection{Result}
In toric geometry, there is a natural way to construct birational maps between toric stacks by presenting toric stacks as GIT quotients (see Section \ref{CTS} for a review) and consider variations of GIT quotients. Certain such VGIT yield crepant birational maps (also known as K-equivalence).

More precisely, let $X_{\pm}$ be toric Deligne-Mumford stacks with $\varphi:X_+ \dashrightarrow X_-$ the birational map arising as the stability parameter crosses a wall. Under certain assumption (see Section \ref{sec:wall_chamber}), there is a common toric blow-up $\Tilde{X}$ with projective birational morphisms $f_{\pm}:\Tilde{X}\ra X_{\pm}$ such that $f_-=\varphi\circ f_+$ and $f^\star_+ K_{X_+}=f^\star_- K_{X_-}$, i.e. $X_+$ and $X_-$ are $K$-equivalent. In \cite{Coates2018}, Coates, Iritani and Jiang considered this situation and proved a crepant transformation correspondence for genus $0$ descendant Gromov-Witten theory of $X_\pm$. Moreover, they relate this correspondence in Gromov-Witten theory with Fourier-Mukai transform between equivariant derived categories of $X_\pm$ studied in \cite{Coates2015a}.

This paper establishes a correspondence in genus $0$ Gromov-Witten theory for crepant birational maps of toric stack {\em bundles} arising from toric wall-crossings. As described in Section \ref{sec:wall_crossing_tsb} below, we consider toric stack bundles $\p_{\pm}$ with toric fiber $X_{\pm}$ related by a crepant toric wall-crossing. We have toric stack bundles $\Tilde{\p}$ with fiber $\Tilde{X}$ and $\p_0$ with fiber $X_0$, where $\Tilde{X}$ is the common toric blow-up of $X_{\pm}$ and $X_0$ is the common blow-down. We have the following diagram with map $\varphi$ induced from the toric wall-crossing on the fiber, $f_{\pm}$ and $g_{\pm}$ induced from the toric blow-ups and blow-downs on the fiber:
\begin{center}
    \begin{tikzcd}
     & \tilde{\p}\arrow[dl,"f_+"']\arrow[dr,"f_-"] &  \\
    \p_+ \arrow[rr,dotted,"\varphi"]\arrow[dr,"g_+"'] & & \p_-\arrow[dl,"g_-"]\\
    & \p_0. &
    \end{tikzcd}
\end{center}

The main result of this paper is the following:

\begin{theorem}\label{thm:main}
    Let $\Tilde{\h}(\p_{\pm})$ be multi-valued version of Givental's symplectic vector space for $\p_{\pm}$. There exist a degree-preserving $R_T(\!(z^{-1})\!)$-linear symplectic transformation $\U:\Tilde{\h}(\p_-)\ra \Tilde{\h}(\p_+)$ such that:
    \begin{enumerate}
        \item $\U$ identifies the $I$-functions of $\p_\pm$: $I_+(y,z)=\U I_-(y,z)$ after analytic continuation;
        \item $\U\circ (g^\star _- v\cup) =(g^\star_+ v\cup) \circ \U$ for all $v\in H^2_T(\p_0)$;
        \item there exists a Fourier-Mukai transformation $\fm:K^0_T(\p_-)\ra K^0_T(\p_+)$ such that the following diagram commutes:
        \begin{equation}\label{diagram:FM_compatible}
    \begin{tikzcd}
    K^0_T(\p_-) \arrow[r,"\fm"]\arrow[d,"{\Psi}_-"] & K^0_T(\p_+) \arrow[d,"{\Psi}_+"] \\
    \Tilde{\h}(\p_-) \arrow[r,"\U"] & \Tilde{\h}(\p_+). 
    \end{tikzcd}
     \end{equation}   
    \end{enumerate}
 Here ${\Psi_\pm=\tilde{\Psi}_\pm\circ \tilde{\tx{ch}}}$ and $\tilde{\Psi}_\pm$ are defined in (\ref{eqn:psi_map}) below.   
\end{theorem}

By \cite{Jiang2017}, the $I$-functions $I_\pm$ of $\p_\pm$, whose definition is recalled in (\ref{eqn:I-function}), lie on the Lagrangian cones of genus $0$ Gromov-Witten theory of $\p_\pm$. Therefore, Theorem \ref{thm:main} can be viewed as a crepant transformation correspondence for genus $0$ Gromov-Witten theory of $\p_\pm$.

Our approach to Theorem \ref{thm:main} is an adaptation of the method used in \cite{Coates2018} for the toric case, replacing the toric mirror theorem of \cite{Coates2015} by its toric stack bundles counterpart in \cite{Jiang2017}. In fact, \cite{Coates2018} is our main motivation for considering Theorem \ref{thm:main} and has been the main source of guidance in our study.

The rest of this paper is organized as follows. In Section \ref{sec:GW_basics} we briefly introduce  Gromov-Witten theory for Deligne-Mumdford stacks and Givental's symplectic formalism. In Section \ref{sec:toric_basics} we recall the definitions and basic results for toric stacks and toric stack bundles, and discuss their toric wall-crossings. In Section \ref{sec:proof_main} we prove Theorem \ref{thm:main}.
\subsection{Acknowledgment}
We thank the referee for helpful comments and suggestions. 
H.-H. T. is supported in part by Simons foundation collaboration grant.

\section{Preliminaries on Gromov-Witten Invariants}\label{sec:GW_basics}

This Section is a collection of basic facts about Gromov-Witten theory. The purpose of this Section is mainly to set up notations we need. Our presentation closely follows those in related papers including \cite{Coates2018}.

\subsection{Chen-Ruan Orbifold Cohomology}

Let $X$ be a smooth proper Deligne-Mumford stack of finite type over $\C$. Let $\abs{X}$ be its coarse moduli space and $IX= X\times_{\abs{X}} X$ the inertia stack of $X$. 
 The stack $IX$ admits a decomposition into connected components
$$IX=\bigsqcup_{v\in \mathsf{B}} X_v,$$
with a ``main component'' $X_0\simeq X$. 
There is a involution map $\tx{inv}:IX\ra IX$ 
which restricts to an isomorphism between corresponding components $\tx{inv}_i: X_i \ra X_{\tx{inv}(i)}$, and becomes identity on $X_0$.

The Chen-Ruan orbifold cohomology $H^{\bullet}_{\tx{CR}}(X)$ of $X$ is the cohomology of the inertia stack $IX$ with grading shifted by a rational number $\iota_v$ called age:
$$H^{k}_{\tx{CR}}(X)=\bigoplus_{v\in \mathsf{B}: k-2\iota_v\in 2\Z} H^{k-2\iota_v}(X_v;\C)$$
$H^{\bullet}_{\tx{CR}}(X)$ carries the orbifold Poincar\'e pairing: for $\alpha,\beta\in H^{\bullet}_{\tx{CR}}(X)$, define
$$(\alpha,\beta)_{\tx{CR}}=\int_{IX} \alpha\cup \tx{inv}^* \beta.$$

Assume $\alpha\in H^\bullet(X_i)$ and $\alpha\in H^\bullet(X_j)$, then $(\alpha,\beta)_{\tx{CR}}$ is only non-zero when $j= \tx{inv}(i)$, in which case
$$(\alpha,\beta)_{\tx{CR}}=\int_{X_i} \alpha\cup \tx{inv}^* \beta$$

If $X$ admits\footnote{When $X$ admits a $T$-action, we can relax the properness requirement: we only assume that the $T$-fixed locus $X^T\subset X$ is proper.} an action of an algebraic torus $T\cong (\C^\times)^m$, we can define the $T$-equivariant Chen-Ruan orbifold cohomology of $X$, $H^{\bullet}_{\tx{CR},T}(X)$, to be the $T$-equivariant cohomology of the inertia stack with age shifting:
$$H^{k}_{\tx{CR},T}(X)=\bigoplus_{v\in \mathsf{B}: k-2\iota_v\in 2\Z} H^{k-2\iota_v}_{T}(X_v;\C).$$
$H^{\bullet}_{\tx{CR},T}(X)$ carries the $T$-equivariant orbifold Poincar\'e pairing: for $\alpha,\beta\in H^{\bullet}_{\tx{CR},T}(X)$, define
$$(\alpha,\beta)_{\tx{CR}}=\int_{IX} \alpha\cup \tx{inv}^* \beta\in S_T$$
using the Atiyah-Bott localization formula. Here $S_T$ is the localization of $R_T=H^{\bullet}_T (\tx{pt};\C)$ with respect to the set of non-zero homogeneous elements.

\subsection{Gromov-Witten Invariant and Quantum Cohomology}
Let $X_{g,n,d}$ be the moduli stack of $n$-pointed genus-$g$ degree-$d$ orbifold stable maps to $X$, where $d\in H_2(\abs{X};\Z)$. Let $\tx{ev}_i:X_{g,n,d}\ra IX$ be the evaluation map at the $i$th marked point, and $\psi_i\in H^2_T(X_{g,n,d})$ the equivariant first Chern class of the $i$-th universal cotangent line bundle over $X_{g,n,d}$. When $X$ admits a $T$-action, $X_{g,n,d}$ carries an induced $T$-action. Consider the $T$-equivariant weighted virtual fundamental class $[X_{g,n,d}]^{\tx{vir}}\in H_{\bullet,T}(X_{g,n,d};\Q)$, see e.g. \cite{Abramovich2008} for details. 

Given $\alpha_1,\dots,\alpha_n \in H^{\bullet}_{\tx{CR},T}(X)$ and non-negative integers $k_1,\dots,k_n$, we have the associated $T$-equivariant Gromov-Witten invariant, defined\footnote{If $X$ is not proper but $X^T$ is proper, this is defined using localization formula for virtual fundamental classes.} to be
$$\langle \alpha_1 \psi^{k_1},\dots,\alpha_n \psi^{k_n}\rangle^{X}_{g,n,d}= \int_{[X_{g,n,d}]^{\tx{vir}}} \prod_{i=1}^n (\tx{ev}^*_i \alpha_i) \psi_i^{k_i}\in S_T.$$

We introduce generating functions of Gromov-Witten invariants. Let $\tx{NE}(X)$ be the cone in $H_2(\abs{X};\R)$ generated by effective curves classes and set 
$$\tx{NE}(X)_{\Z}=\{d\in H_2(\abs{X};\Z):d\in \tx{NE}(X)\}.$$ 
Denote by $Q$ the Novikov variable, so for a ring $R$, we have $$R\llbracket Q \rrbracket=\left\{\sum_{d\in \tx{NE}(X)_\Z} a_d Q^d : a_d\in R\right\}.$$
For $\mathbf{t}=\mathbf{t}(z)=t_0+t_1 z+t_2 z^2+\cdots \in H^{\bullet}_{\tx{CR},T}(X)\otimes_{R_T} S_T[z]\llbracket Q \rrbracket$, define
$$\langle \mathbf{t},\dots,\mathbf{t}\rangle^X_{g,n,d} = \langle \mathbf{t}(\psi),\dots,\mathbf{t}(\psi)\rangle^X_{g,n,d}= \sum_{k_1,\dots,k_n\geq 0} \langle t_{k_1} \psi^{k_1},\dots,t_{k_n} \psi^{k_n}\rangle^{X}_{g,n,d}.$$
The genus-$g$ descendant potential of $X$ is defined to be
$$\mathcal{F}^g_X(\mathbf{t})=\sum_{n\geq 0} \sum_{d\in \tx{NE}(X)_\Z } \frac{Q^d}{n!} \langle \mathbf{t},\dots,\mathbf{t}\rangle^X_{g,n,d}.$$

\subsection{Equivariant Quantum Cohomology}
 Let $\{\phi_i\}$ be a basis of $H_{\tx{CR}}^\bullet(X)$. Let $\eta_{ij}=(\phi_i,\phi_j)_{\tx{CR}}$ and $\eta^{ij}$ be the entries of inverse matrix of $(\eta_{ij})$. 
For $\tau\in H_{\tx{CR}}^\bullet(X)$, we define the quantum product at $\tau$ by
\begin{align*}
    \phi_i\bullet_\tau \phi_j &= \sum_{k,l} \sum_{d\in \tx{NE}(X)_\Z} \sum_{n=0}^\infty \frac{Q^d}{n!}\langle \phi_i,\phi_j,\phi_k,\tau,\dots,\tau \rangle_{0,n+3,d}^X \eta^{kl}\phi_l.
\end{align*}
Equivalently, 
\begin{align*}
    (\phi_i\bullet_\tau \phi_j,\phi_k)= \sum_{d\in \tx{NE}(X)_\Z} \sum_{n=0}^\infty \frac{Q^d}{n!}\langle \phi_i,\phi_j,\phi_k,\tau,\dots,\tau \rangle_{0,n+3,d}^X.
\end{align*}

For equivariant orbifold cohomology $H^\bullet_{\tx{CR},T}(X)$ and $\tau\in H^\bullet_{\tx{CR},T}(X)$, we define the quantum product analogously, so
\begin{align*}
    (\phi_i\bullet_\tau \phi_j,\phi_k)_{\tx{CR}}= \sum_{d\in \tx{NE}(X)_\Z} \sum_{n=0}^\infty \frac{Q^d}{n!}\langle \phi_i,\phi_j,\phi_k,\tau,\dots,\tau \rangle_{0,n+3,d}^X,
\end{align*}
or
\begin{align*}
    \phi_i\bullet_\tau \phi_j 
    &= \sum_{d\in \tx{NE}(X)_\Z} \sum_{n=0}^\infty \frac{Q^d}{n!} \tx{inv}^* \tx{ev}_{3,*} \left( \tx{ev}^*_1(\phi_i) \tx{ev}^*_2(\phi_j) \prod_{l=4}^{n+3} \tx{ev}^*_l(\tau) \cap [X_{0,n+3,d}]^{\tx{vir}} \right).
\end{align*}

The quantum product $\bullet_\tau$ defines an associative and commutative ring structure on $H^{\bullet}_{\tx{CR},T}(X)\otimes_{R_T} R_T[z]\llbracket \tau, Q \rrbracket$, which is called the equivariant quantum cohomology of $X$.

\subsection{Givental's Symplectic Formalism}
The original reference for Givental's formalism is \cite{Givental2004SymplecticStructures}. See \cite{Coates2009ComputingInvariants} for a presentation of its orbifold version.

Givental's symplectic space is defined to be
$$\h=H^{\bullet}_{\tx{CR},T}(X)\otimes_{R_T} S_T(\!(z^{-1})\!)\llbracket Q \rrbracket,$$
where $S_T(\!(z^{-1})\!)$ is the ring of formal Laurent series in $z^{-1}$ with coefficients in $S_T$. The symplectic form on the space is given by 
$$\Omega(f,g)=-\tx{Res}_{z=\infty} (f(-z),g(z))_{\tx{CR}} dz,$$
for any $f,g\in \h$, i.e. the coefficient of $z^{-1}$ of the orbifold Poincar\'e pairing of $f(-z)$ and $g(z)$. This space has a standard polarization
$$\h=\h_+ \oplus \h_-,$$
where 
$$\h_+=H^{\bullet}_{\tx{CR},T}(X)\otimes_{R_T} S_T[z]\llbracket Q \rrbracket,$$
$$\h_-=z^{-1}H^{\bullet}_{\tx{CR},T}(X)\otimes_{R_T} S_T\llbracket z^{-1} \rrbracket\llbracket Q \rrbracket.$$
The subspaces $\h_{\pm}\subset \h$ are the maximal isotropic subspaces with respect to $\Omega$, and $\Omega$ induces a non-degenerate pairing between $\h_+$ and $\h_-$. This identifies $\h=\h_+\oplus\h_-$ with the total space of the cotangent bundle $T^*\h_+$ of $\h_+$.

Let $\{\phi^i\},\{\phi_i\}\subset H^{\bullet}_{\tx{CR},T}(X)\otimes_{R_T} S_T$ be bases dual with respect to orbifold Poincar\'e pairing, i.e. $(\phi^i,\phi_j)_{\tx{CR}}=\delta^i_j$. A general point in $\h$ takes the form
$$\sum_{a=0}^\infty \sum_{i} p_{a,i} \phi^i (-z)^{-a-1} + 
\sum_{b=0}^\infty \sum_{j} q_{b}^{j} \phi_j z^b.$$
This defines Darboux coordinates $\{p_{a,i},q_b^j\}$ on $(\h,\Omega)$ compatible with the polarization $\h=\h_+\oplus \h_-$. Let $p_a=\sum_{i} p_{a,i} \phi^i, q_b=\sum_{j} q_b^j \phi_j$, and denote 
$$\mathbf{p}=\sum_{k=0}^\infty p_k (-z)^{-k-1}, \quad 
 \mathbf{q}=\sum_{k=0}^\infty q_k z^k.$$
Givental's Lagrangian cone $\mathcal{L}_X$ is defined to be the graph of the differential of the genus-zero descendant potential $\mathcal{F}^0_X$,
$$\mathcal{L}_X=\left\{ (\mathbf{p},\mathbf{q})\in \h_+\oplus \h_- : \mathbf{p}= d_{\mathbf{q}} \mathcal{F}^0_X \right\}\subset T^*\h_+ \cong \h.$$
Let $x=(x_1,\dots,x_m)$ be formal variables, an $S_T \llbracket Q, x \rrbracket$-valued point of $\mathcal{L}_X$ is an element of the form
$$-1z+\mathbf{t}+\sum_{n=0}^\infty \sum_{d\in \tx{NE}(X)_\Z } \sum_{i} \frac{Q^d}{n!} \left\langle \mathbf{t},\dots,\mathbf{t}, \frac{\phi_i}{-z-\psi}\right\rangle^X_{0,n+1,d} \phi^i$$
for some $\mathbf{t}=\mathbf{t}(z)\in \h_+\llbracket x \rrbracket$ satisfying
$$\mathbf{t}|_{x=Q=0}=0.$$

Let $J$-function be the following formal series
$$J_X(\mathbf{t},z)= 1z+\mathbf{t}+\sum_{n=0}^\infty \sum_{d\in \tx{NE}(X)_\Z } \sum_{i} \frac{Q^d}{n!} \left\langle \mathbf{t},\dots,\mathbf{t}, \frac{\phi_i}{z-\psi}\right\rangle^X_{0,n+1,d} \phi^i.$$
Then $J_X(\mathbf{t},-z)$ gives an $S_T \llbracket Q, x \rrbracket$-valued point of the Lagrangian cone $\mathcal{L}_X$.

Givental's Lagrangian cone $\mathcal{L}_X$ has an important geometric property: each tangent space $\mathcal{T}$ to the cone is tangent to the cone exactly along $z\mathcal{T}$, see \cite{Givental2004SymplecticStructures}. Using the special geometric properties of $\mathcal{L}_X$, we can reconstruct the quantum product $\bullet_\tau$ from the Lagrangian cone $\mathcal{L}_X$ and vice versa. 

For instance, let $v\in H^\bullet_{\tx{CR},T}(X)$, the equivariant quantum connection
$$\nabla_v:H^{\bullet}_{\tx{CR},T}(X)\otimes_{R_T} R_T[z]\llbracket \tau, Q \rrbracket \ra z^{-1}H^{\bullet}_{\tx{CR},T}(X)\otimes_{R_T} R_T[z]\llbracket \tau, Q \rrbracket$$
is defined by 
$$\nabla_v f(\tau) = \partial_v f(\tau) + z^{-1} v\bullet_\tau f(\tau).$$
Denote $\nabla_{\phi_i}$ by $\nabla_i$. Then the fundamental solution 
$$L(\tau,z)\in \tx{End}_{R_T}(H^{\bullet}_{\tx{CR},T}(X))\otimes_{R_T} R_T (\!(z^{-1})\!) \llbracket \tau, Q \rrbracket$$
is determined by the following conditions
\begin{align*}
    \nabla_i L(\tau,z)\phi &= 0,\\
    \left( vQ\frac{\partial}{\partial Q}-\partial_v \right)L(\tau,z)\phi &= L(\tau,z)\frac{v}{z} \phi,\\
    L(\tau,z)|_{\tau=Q=0} &= \tx{id},
\end{align*}
for $i$ runs over the basis $\{\phi_i\}$ of $H^\bullet_{\tx{CR},T}(X)$, $\phi\in H^\bullet_{\tx{CR},T}(X)$, and $v\in H^2_T(X)$.

The fundamental solution can be written in the following form using Gromov-Witten invariants:
$$L(\tau,z)\phi_i = \phi_i+ \sum_{j} \sum_{d\in \tx{NE}(X)_\Z } \sum_{n=0}^\infty   \frac{Q^d}{n!} \left\langle \frac{\phi_i}{-z-\psi}, \tau ,\dots,\tau, \phi_j \right\rangle^X_{0,n+2,d} \phi^j.$$
The map $\tau\mapsto T_\tau = L(\tau,-z)^{-1} \h_+$ gives a versal family of tangent spaces to Givental's Lagrangian cone $\mathcal{L}_X$, and the cone can be reconstructed by $\mathcal{L}_X = \bigcup_\tau z T_\tau$.

\section{Preliminaries on Toric Geometry}\label{sec:toric_basics}

In this Section we present constructions and basic results on toric stacks and toric stack bundles. Our presentation closely follows \cite{Coates2018}.

\subsection{Toric Variety Associated To A Fan}
Fix lattices $N, M$ dual to each other and denote $N_{\R}=N\otimes \R, M_{\R}= M\otimes \R$. A convex polyhedral cone, or simply a cone, is a set in $N_{\R}$ of the form
$$\sigma = \tx{Cone}(S) = \left\{ \sum_{u\in S} \lambda_u u: \lambda_u\geq 0   \right\},$$
where $S\subset N_{\R}$ is finite and call $\sigma$ a cone generated by $S$. The dual cone of $\sigma$ is defined to be
$$\sigma^\vee =\left\{  m\in M_{\R}: \langle m,u\rangle \geq 0 \ \tx{for all}\ u\in \sigma   \right\},$$
which is a cone in $M_{\R}$ and $(\sigma^\vee)^\vee = \sigma$. We define the hyperplane 
$$H_m= \left\{ u\in N_{\R}: \langle m, u \rangle=0 \right\}.$$
A face of the cone $\sigma$ is $\tau = H_m \cap \sigma$ for some $m\in \sigma^\vee$, denoted by $\tau \preceq \sigma$. Every face of $\sigma$ is again a cone. Moreover, an intersection of faces or face of a face is again a face of $\sigma$. A face of codimension $1$ is called a facet and a face of dimension $1$ is called an edge.

A cone is called rational if the generating set $S$ is in $N\subset N_\mathbb{R}$. A strongly convex rational cone has minimal generators called the ray generators of the edges. We call such a cone simplicial if its minimal generators are linearly independent over $\R$, and call such a cone smooth if the minimal generators form part of a $\Z$-basis of $N$.

A rational cone $\sigma\subset N_{\R}$ gives a semigroup $S_{\sigma}=\sigma^\vee \cap M$, and the affine variety associated to $\sigma$ is $$U_{\sigma} = \tx{Spec}(\C[S_{\sigma}]).$$

A fan $\Sigma$ in $N_{\R}$ is a finite collection of rational cones $\sigma\subset N_{\R}$ such that for each $\sigma\in \Sigma$, each face of $\sigma$ is also in $\Sigma$, and for any $\sigma_1, \sigma_2 \in \Sigma$, the intersection $\sigma_1\cap \sigma_2$ is a face of each, hence also in $\Sigma$. The support of the fan is $|\Sigma|= \bigcup_{\sigma\in\Sigma} \sigma \subset N_{\R}$.

Given a fan $\Sigma$ in $N_{\R}$, the associated toric variety $X_{\Sigma}$ is defined by gluing the affine toric varieties $U_{\sigma}, U_{\sigma'}$ along $\sigma\cap \sigma'$ for any $\sigma,\sigma'\in \Sigma$. In fact, any normal separated toric variety comes from a fan $\Sigma$. 

Geometric information of a toric variety can be deduced from the combinatorial data of its fan. For example, $X_{\Sigma}$ is smooth if and only if every cone $\sigma \in \Sigma$ is smooth; $X_{\Sigma}$ is an orbifold if and only if every cone $\sigma \in \Sigma$ is simplicial; $X_{\Sigma}$ is compact if and only if $\Sigma$ is complete, i.e. the support $|\Sigma|= \bigcup_{\sigma\in\Sigma} \sigma $ is the entire $ N_{\R}$.

\subsubsection{Torus Fixed Points and Toric Blow-ups}
There is an orbit-cone correspondence for toric varieties. Let $X_{\Sigma}$ be the toric variety associated to a fan $\Sigma$ in $N_{\R}$. $X_{\Sigma}$ contains the torus $T_N = N\otimes \C^\times$ as an open subset and $X_\Sigma$ has a $T_N$-action which extends the action of $T_N$ on itself. 

Let $\sigma$ be a cone in $\Sigma$, $N_\sigma$ the sublattice of $N$ spanned by the points in $\sigma\cap \N$ and $N(\sigma)=N/N_{\sigma}$. The cones $\sigma \in \Sigma$ are in bijective correspondence with $T_{N}$-orbits $O(\sigma)$ of $X_{\Sigma}$. Recall the definition of $U_\sigma$, a point $p\in U_\sigma$ can be represented by a map $\gamma: S_\sigma \ra \C$ which sends $m\in S_\sigma$ to $\chi^m (p) \in \C$. The torus orbit is
\begin{align*}
    O(\sigma) &= \left\{ \gamma: S_\sigma \ra \C | \gamma(m)\neq 0 \Leftrightarrow m\in \sigma^{\perp}\cap M \right\}\simeq \tx{Hom}(\sigma^{\perp}\cap M, \C^\times) \simeq T_{N(\sigma)}.
\end{align*}
In particular, a top-dimensional cone corresponds to a torus-fixed point of $X_{\Sigma}$.

Operations on a fan $\Sigma$ yield morphisms for the associated toric variety $X_\Sigma$. Given a fan $\Sigma$ in $N_{\R}$, a fan $\Sigma'$ refines $\Sigma$ if every cone of $\Sigma'$ is contained in a cone of $\Sigma$ and $|\Sigma'|=|\Sigma|$. Moreover, let $\sigma$ be a smooth top-dimensional cone in $\Sigma$ generated by $u_1,\dots,u_n$, then $u_1,\dots,u_n$ form a basis for $N$. Let $u_0=u_1+\cdots+u_n$ and $\Sigma'(\sigma)$ the set of cones generated by subsets of $\{u_0,u_1,\dots,u_n\}$ not containing $\{u_1,\dots,u_n\}$. Define
$$\Sigma^*(\sigma) = (\Sigma\setminus \{\sigma\}) \cup \Sigma'(\sigma),$$
the star subdivision of $\Sigma$ along $\sigma$. This is again a fan in $N_{\R}$. 

The star subdivision $\Sigma^*(\sigma)$ is a refinement of $\Sigma$, and the induced toric morphism 
$$\phi: X_{\Sigma^*(\sigma)} \ra X_{\Sigma}$$
is the blow-up of $X_{\Sigma}$ at the distinguished point corresponding to $\sigma$.

Similarly, we can define a toric blow-up along a subvariety $V(\tau)$ using the star subdivision of $\Sigma$ along $\tau$. Let $\tau \in \Sigma$ be such that all cones of $\Sigma$ containing $\tau $ are smooth. Let $u_{\tau} = \sum_{\rho\in \tau(1)} u_{\rho}$ and for each $\sigma\in \Sigma$ containing $\tau$, define
$$\Sigma^*_{\sigma}(\tau) = \left\{ \tx{Cone}(A): A\subseteq \{u_{\tau} \}\cup \sigma(1), \tau(1)\nsubseteq A\right\},$$
and define the star subdivision of $\Sigma$ along $\tau$ to be
$$\Sigma^*(\tau) = \left\{ \sigma\in \Sigma: \tau \nsubseteq \sigma  \right\} \cup \bigcup_{\tau\subseteq \sigma}  \Sigma^*_{\sigma}(\tau).$$
Then the induced toric morphism 
$$\phi: X_{\Sigma^*(\tau)} \ra X_{\Sigma}$$
is the blow-up of $X_{\Sigma}$ along  $V(\tau)$.

\subsection{Toric Stacks}\label{CTS}
Given $K\cong (\C^\times)^r$, characters of $K$, $$D_1,\dots,D_m \in \Ll^\vee = \text{Hom}(K,\C^\times),$$ define a map from $K$ to $(\C^\times)^m$, hence an action of $K$ on $\C^m$. For a set $I\subset [m]=\{1,2,\dots,m\}$, define
\begin{align*}
\angle_{I} &=\left\{\sum_{i\in I} a_i D_i : a_i\in \R, a_i>0\right\} \subset \Ll^\vee \otimes \R,\\
(\C^\times)^I \times (\C)^{\ol{I}} &= \left\{(z_1,\dots,z_m): z_i\neq 0 \ \tx{for} \ i\in I \right\} \subset \C^m,
\end{align*}
where $\ol{I}=[m]\setminus I$ is the complement of $I$. Set $\angle_{\emptyset}=\{0\}$.

Given a stability condition $\omega \in \Ll^\vee \otimes \R$, we define\footnote{Elements of $\A_\omega$ are often called {\em anticones}.}
\begin{align*}
\A_\omega = \{ I\subset [m]: \omega\in \angle_I \},\quad 
U_\omega = \bigcup_{I\in \A_\omega} (\C^\times)^I \times \C^{\ol{I}},\quad 
X_\omega = [U_\omega/K],
\end{align*}
with the assumptions that $[m]\in \A_\omega$ and for each $I\in \A_\omega$, $\{ D_i: i\in I\}$ spans $\Ll^\vee \otimes \R$ over $\R$. 

Let $S\subset [m]$ be the set of indices $i$ such that $[m]\setminus \{i\} \notin \A_\omega$. Then the characters $\{D_i: i\in S\}$ are linearly independent and every $I\in \A_\omega$ contains $S$. In fact we may define $X_\omega=[U'_\omega/G]$ where $\A_\omega =\{I\sqcup S: I\in \A'_\omega\}$, $U_\omega\cong U'_\omega \times (\C^\times)^{|S|}$, $G=\tx{Ker}(K\ra (\C^\times)^{|S|}) $ and get the same $X_\omega$. And this is the original construction of toric Deligne-Mumford stacks of Borisov-Chen-Smith \cite{Borisov2004TheStacks}. 

The definition here is actually equivalent to the construction by Jiang using an $S$-extended stacky fan \cite{Yunfeng2008}. An $S$-extended stacky fan $\mathbf{\Sigma}=(\mathbf{N},\Sigma,\beta,S)$ consists of the following data: a finitely generated abelian group $\mathbf{N}$, a rational simplicial fan $\Sigma$ in $\mathbf{N}\otimes \R$, a homomorphism $\beta: \Z^m\ra \N$ with $b_i=\beta(e_i)\in \N$ and $\ol{b_i}$ the image of $b_i$ in $\N\otimes\R$, and $S\subset [m]$ such that $1$-dimensional cones in $\Sigma$ are in  bijective correspondence with $\ol{b_i}$ for $i\in [m]\setminus S$, and for each $i\in S$, $\ol{b_i}$ lies in the support $\abs{\Sigma}$.

We now explain the equivalence. Starting from the GIT data, consider the following exact sequence
$$0\ra \Ll \overset{(D_1,\dots,D_m)}{\xrightarrow{\hspace{1cm}}} \Z^m \overset{\beta}{\longrightarrow} \N \ra 0, $$
i.e. define $\beta:\Z^m \ra \N$ to be the cokernel of $\Ll\ra \Z^m$. Let $b_i=\beta(e_i)$ and define $S$ as the set of indices $i$ such that $[m]\setminus \{i\} \notin \A_\omega$. Then $\Sigma_\omega$ is the fan generated by elements $\{\ol{b_i}, i\in [m]\setminus S\}$ and we get the $S$-extended stacky fan $\mathbf{\Sigma}=(\mathbf{N},\Sigma_\omega,\beta,S)$.

Conversely, starting from an $S$-extended stacky fan. Using the same exact sequence above, we may define $\Ll$ as the kernel of $\beta$ and $K=\Ll\otimes \C^\times$. Take the dual exact sequence 
$$0\ra \N^\vee \longrightarrow  (\Z^m)^\vee \longrightarrow \Ll^\vee,$$
and define $D_i\in \Ll^\vee$ to be the image of $e_i^\vee$ in the dual sequence. Define
$$\A_\omega=\{I\subset[m]: S\subset I \ \tx{and}\ \sigma_{\ol{I}} \ \tx{is a cone of} \  \Sigma\},$$
and take $\omega$ in $\cap_{I\in \A_\omega} \angle_I$, then this recovers the GIT data.

\subsubsection{Equivariant Cohomology of Toric Stacks}\label{ECTS}
Here we summarize some results about equivariant cohomology of toric stacks.

Let $\q=T/K$, then the $T$-action on $U_\omega$ descends to a $\q$-action on $X_\omega=[U_\omega/K]$. Let $R_\q=H^{\bullet}_\q(\tx{pt};\C)$, then $R_\q\cong \tx{Sym}^\bullet (\mathbf{N}^\vee \otimes \C)$. Let $\lambda_i$ be the equivariant first Chern class of the irreducible $T$-representation given by the projection $T\cong (\C^\times)^m \ra \C^\times$ to the $i$th factor, then $R_T=\C[\lambda_1,\dots,\lambda_m]$.

The $\q$-equivariant cohomology ring of $X_\omega$ can be written as
$$H^\bullet_\q (X_\omega;\C) = R_\q [u_1,\dots,u_m]/(\mathfrak{I}+\mathfrak{J}),$$
where $u_i$ is the $\q$-equivariant class Poincar\'e dual to the toric divisor 
$$[\left\{ (z_1,\dots,z_m)\in U_\omega: z_i=0 \right\}/K],$$
and $\mathfrak{I},\mathfrak{J}$ are ideals defined by 
$$\mathfrak{I}=\left \langle \chi-\sum_{i=1}^m\langle \chi,b_i\rangle u_i: \chi\in \mathbf{N}^\vee\otimes \C \right\rangle,\quad 
\mathfrak{J}=\left \langle \prod_{i\notin I} u_i: I\notin \A_\omega \right \rangle.$$
Note the relations in $\mathfrak{J}$ imply that $u_i=0$ in $H^\bullet_\q (X_\omega;\C)$ for $i\in S$.

The $T$-equivariant cohomology ring of $X_\omega$ is given by the extension of scalars, 
$$H^\bullet_T (X_\omega;\C)\cong H^\bullet_\q (X_\omega;\C)\otimes_{R_\q} R_T,$$
where the homomorphism $R_\q\ra R_T$ is given by $\chi \mapsto \sum_{i=1}^m\langle \chi,b_i\rangle \lambda_i$ for $\chi\in \mathbf{N}^\vee \otimes \C$. We also regard $u_i$ as a $T$-equivariant class under this identification.

There is a split short exact sequence
\begin{equation}\label{eqn:split_exact_seq}
0\ra H^2_\q(X_\omega;\R)\ra H^2_T(X_\omega;\R)\ra H^2_K(\tx{pt};\R)\ra 0,
\end{equation}
where the splitting is given by $\theta:H^2_K(\tx{pt};\R)\cong \Ll^\vee\otimes \R \ra H^2_T(X_\omega;\R)$, $\theta(D_i)=u_i-\lambda_i$. The class $\theta(p)$ actually gives the $T$-equivariant first Chern class of a line bundle $L(p)$ determined by the $K$-representation $p\in \Ll^\vee$.

\subsection{Wall and Chamber Structure}\label{sec:wall_chamber}
The space of stability conditions $\Ll^\vee \otimes \R$ has a wall and chamber structure, which describes how the GIT quotient $X_\omega$ changes as $\omega\in \Ll^\vee \otimes \R$ varies. 

For each $\omega$, the chamber containing $\omega$ is given by 
$$C_\omega = \bigcap_{I\in \A_\omega} \angle_I.$$
The quotient $X_{\omega'}$ is isomorphic to $X_\omega$ for any $\omega'\in C_\omega$. However $X_{\omega'}$ changes when $\omega'$ crosses a codimension one boundary (called a wall) of $C_\omega$.

We also define the cone of ample divisors 
$$C'_\omega = \bigcap_{I\in \A'_\omega} \angle'_I,$$
where $\A'_\omega=\{I\setminus S: I\in \A_\omega\}$, $\angle'_I=\sum_{i\in I} \R_{>0} D'_i$ and $D'_I$ is the image of $D_i$ in $\Ll^\vee\otimes \R/ \sum_{i\in S} \R D_i\cong H^2(X_\omega;\R)$. Then
$$C_\omega\cong C'_\omega \times \left( \sum_{i\in S} \R_{>0} D_i \right).$$

Let $C_+$ and $C_-$ be two chambers in $\Ll^\vee\otimes\R$ separated by a hyperplane $W$. Pick  stability conditions $\omega_+\in C_+,\omega_-\in C_-$. Consider $\A_{\pm}=\A_{\omega_{\pm}}$ as defined above, and the corresponding toric stacks $X_{\pm}=X_{\omega_{\pm}}$. Let $\varphi: X_+\dashrightarrow X_-$ be the birational transformation induced by the toric wall-crossing. We make the additional assumption: 
$$\sum_{i=1}^m D_i \in W,$$
so that transformation $\varphi$ is crepant. Let $e\in W^{\perp}\cap \Ll$ be a primitive lattice vector such that $e$ is positive on $C_+$ and negative in $C_-$.

The chambers in $\Ll^\vee \otimes \R$ form a fan called the Gelfand-Kapranov-Zelevinsky (GKZ) fan. Especially, we consider the subfan consisting of cones $\overline{C_+}$ and $\overline{C_-}$. Let $\M$ be the toric variety associated to this fan. Since $\overline{C_\pm}$ are in general not simplicial cones, $\M$ may be singular. So we may work in some overlattices of $\Ll$. Define 
$$\K_\pm=\left\{ f\in \Ll\otimes\Q:I_f\in \A_{\pm}\right\},$$ where $I_f=\{i\in [m]: D_i \cdot f\in \Z\}$. Let $\Tilde{\Ll}_\pm$ be the free $\Z$-submodule generated by $\K_\pm$. Note that $D_j$ is always in $\tilde{\Ll}_\pm^\vee$ for any $j\in S_{\pm}$ and they are linearly independent, so we may extend $\{D_j:j\in S_{\pm}\}$ to a basis of $\tilde{\Ll}_\pm^\vee$. After reordering, the basis can be written as $\left\{  p^{\pm}_1,\dots,p^{\pm}_{r-1},p^{\pm}_r  \right\}$
such that $p^+_i=p^-_i \in W$ for $1\leq i\leq r-1$ and $p^\pm_r$ is the unique vector in each basis that does not lie on $W$. Let $\{y^\pm_1,\dots,y^\pm_r\}$ be the corresponding coordinates, and, for any $d\in \Ll\otimes \R$, define
$$(y^{\pm})^d= \prod_{i=1}^r (y^{\pm}_i)^{p^\pm_i \cdot d}.$$
Under change of coordinates, we have $(y^+)^e = (y^-)^e$.

\begin{remark}
One can interpret $\mathcal{M}$ as the moduli space of certain Landau-Ginzburg models (mirror to the toric stacks). For this paper, our use of $\mathcal{M}$ is mainly to house $y^\pm$, which are variables in  $I$-functions. See \cite{Coates2018} for a detailed discussion about $\mathcal{M}$.
\end{remark}

\subsection{Toric Stack Bundles}\label{sec:tsb}
 Let $P\ra B$ be a principal $(\C^\times)^m$-bundle over a smooth projective variety $B$. Let $K, D_i$ and $U_\omega$ be as in Section \ref{CTS}. Associated to these data, we define a toric stack bundle to be the quotient
 \begin{equation}\label{eqn:toric_stack_bundle}
 \p = {^P}X_\omega =[\left(P\times_{(\C^\times)^m} U_\omega\right)/K ]\to B,    
 \end{equation} 
where $K$ acts on $U_\omega$ via characters $D_i$ and acts on $P$ trivially. $\mathcal{P}$ is a fiber bundle with fiber $X_w$ over $B$ and it is constructed from the principal $(\mathbb{C}^\times)^m$-bundle $P\to B$ by changing fibers. More details about construction of toric stack bundles can be found in \cite{JT}.

Let $\mathbf{\Sigma}_\omega=(\mathbf{N},\Sigma_\omega,\beta,S)$ be the corresponding $S$-extended stacky fan corresponding to $K, \{D_i\},\omega$ as recalled in Section \ref{CTS}, so $X_\omega=X(\mathbf{\Sigma}_\omega)$. The inertia stack of $\p$ is indexed by elements of the following set
$$\tx{Box}(\mathbf{\Sigma}_\omega)=\left\{ v\in \mathbf{N}: \overline{v}=\sum_{i\notin I} c_i \overline{b_i} \in \mathbf{N}\otimes \R \ \tx{for some} \ I\in\A, 0\leq c_i<1 \right\},$$
via the decomposition 
$$I\p=\coprod_{v\in \tx{Box}(\mathbf{\Sigma}_\omega)} \p_v = \coprod_{v\in \tx{Box}(\mathbf{\Sigma}_\omega)} {^P}X(\mathbf{\Sigma}_\omega/\sigma(\overline{v})). $$
Now define $\K=\left\{ f\in \Ll\otimes\Q:I_f\in \A_{\omega}\right\}$ where $I_f=\{i\in [m]: D_i \cdot f\in \Z\}$. There is a bijection between $\K/\Ll$ and $\tx{Box}(\mathbf{\Sigma}_\omega)$ given by
$$[f] \leftrightarrow v_f= \sum_{i=1}^m \lceil -(D_i\cdot f)\rceil b_i. $$
So we can rewrite the decomposition of $I\p$ over $\K/\Ll$ as
$$I\p=\coprod_{f\in \K/\Ll} \p_f. $$
The $T$-equivariant Chen-Ruan cohomology is 
$$H^{\bullet}_{\tx{CR},T}(\p)= \bigoplus_{f\in \K/\Ll} H^{\bullet - 2\iota_f}_T (\p_f), $$
where $\iota_f = \sum_{i\notin I_f} \langle D_i\cdot f \rangle$ is the age.

Now assume $P$ arises from a direct sum of line bundles $\oplus_{j=1}^m L_j$ by removing the zero section, and let $\Lambda_j= c_1(L_j^*)$. We may take $K_\R=(S^1)^r$ and regard the quotient (\ref{eqn:toric_stack_bundle}) as a fiberwise symplectic reduction for the action of $K_\R$. This defines $r$ tautological bundles over $\p$. Define $-P_1,\dots,-P_r$ to be the $T$-equivariant first Chern classes of the tautological bundles, and $p_1,\dots,p_r$ the restrictions of $P_1,...,P_r$ to the fiber $X_\omega$. Recall that $\lambda_1,\dots,\lambda_m$ are the first Chern classes of the $m$ duals of tautological bundles on $BT$ and $R_T:=H^\bullet_T(\tx{pt};\C)=\C[\lambda_1,\dots,\lambda_m]$. The classes $P_1,\dots,P_r$ generate the cohomology ring $H^{\bullet}_T(\p)$ over $H^{\bullet}(B)\otimes R_T$ just as $p_1,\dots,p_r$ generate $H^{\bullet}_T(X_\omega)$ over $R_T$. This gives a basis of $H^{\bullet}_T(\p)$ over $H^{\bullet}(B)\otimes R_T$ consisting of monomials in $\{P_i\}$.

Recall that we have class $u_j$ defined as dual to the $j$-th toric divisor in $X_\omega$. Thus we have $u_j=\sum_{i=1}^r d_{ij}p_i$ if we write $D_j=(d_{ij})_{i=1}^r$. Define
$$U_j=\sum_{i=1}^r d_{ij}P_i-\Lambda_j-\lambda_j,$$
so $U_j$ is the $T$-equivariant cohomology class dual to the $j$-th toric divisor in $\p$. The splitting of the short exact sequence (\ref{eqn:split_exact_seq}) 
gives rise to a map $\theta: H^2(B;\R)\otimes R_K \ra H^2_T(\p;\R)$ which maps $D_j$ to $U_j-\lambda_j$ and fixes classes in $H^2(B)$. Especially, for a line bundle over $\p$ defined by a character $p$, its $T$-equivariant first Chern class is $\theta(p)$.


\subsubsection{Wall-Crossing for Toric Stack Bundles}\label{sec:wall_crossing_tsb}
By the construction of toric stack bundles (\ref{eqn:toric_stack_bundle}), as the stability parameter $\omega$ crosses a wall, $\mathcal{P}$ changes as the fiber $X_\omega$ changes.

As in Section \ref{sec:wall_chamber}, consider chambers $C_{\pm}\subset \Ll\otimes \R$ separated by a wall $W$. Let $\omega_{\pm}\in C_{\pm}$ be stability conditions and let 
$$\p_\pm =[P\times_{(\C^\times)^m} U_{\omega_{\pm}}/K ]$$ 
be the corresponding toric stack bundles. Also choose a stability condition $w_0\in W$ such that it is in the relative interior of $W\cap \overline{C_+}= W\cap \overline{C_-}$. Let $\Tilde{\p}$ be the toric stack bundle with fiber $\Tilde{X}$ which is the common blow-up of $X_{\pm}$, and let $\p_0$ be the toric stack bundle with fiber\footnote{$X_{\omega_0}$ is an Artin stack.} $X_{\omega_0}$. Note that $\p_0$ contains both $\p_{\pm}$ as open substacks. We have the following commutative diagram
\begin{center}
    \begin{tikzcd}
     & \tilde{\p}\arrow[dl,"f_+"']\arrow[dr,"f_-"] &  \\
    \p_+ \arrow[rr,dotted,"\varphi"]\arrow[dr,"g_+"'] & & \p_-\arrow[dl,"g_-"]\\
    & \p_0, &
    \end{tikzcd}
\end{center}
where each map is induced from a map on toric fibers. Especially, the canonical bundles $K_{\p_{\pm}}$ of $\p_\pm$ are restrictions of the same line bundle $L_0$ on $\p_0$ given by a line bundle on $B$ and the character $-\sum_{i=1}^m D_i$. The condition $\sum_{i=1}^m D_i\in W$ ensures that $L_0$ comes from a $\Q$-Cartier divisor $\mathcal{D}_0$ on $\p_0$. Then we have 
$$f^*_+(K_{\p_+})=f^*_+ \circ g^*_+ (\mathcal{D}_0) = f^*_- \circ g^*_- (\mathcal{D}_0) = f^*_-(K_{\p_-}),$$
which shows that the birational map $\varphi$ is crepant.

\section{Crepant Transformation Correspondence}\label{sec:proof_main}

In this Section we prove Theorem \ref{thm:main}. First, we compute the analytic continuation of the $I$-function using the Mellin-Barnes method, which gives the linear transformation $\U$. Then we compute the Fourier-Mukai transformation on $K$-groups and compare it with $\U$.

\subsection{$I$-Function and $H$-Function}
The $I$-function for a toric stack bundle $\p$ is defined in \cite{Jiang2017} to be
\begin{equation}\label{eqn:I-function}
I^{\p}(y,z) =ze^{\sigma/z} \sum_{\D\in \overline{NE}(B)} J^B_{\D}(z,\tau) Q^{\D} \sum_{d\in\K} y^d \Delta \mathbf{1}_{\left [ -d \right ]},
\end{equation}
where $J^B_\D(z,\tau)$ is the degree-$\D$ component of the $J$-function of $B$, $\sigma\in H^2_T(\p)$ is a class given by $\sigma= \sum_{i=1}^r \theta(p_i) \tx{log} (y_i) +c_0(\lambda)$ where $c_0(\lambda)=\lambda_1+\cdots+\lambda_m$ and $\theta$ is the map defined in Section \ref{sec:tsb}, and
\begin{align*}
\Delta &=\prod_{j=1}^m \frac{\prod_{\langle a\rangle =\langle D_j\cdot d+ \Lambda_j \D\rangle, a\leq 0} (U_j+az)}{\prod_{\langle a\rangle =\langle D_j\cdot d+\Lambda_j \D\rangle, a\leq D_j\cdot d+\Lambda_j \D} (U_j+az)}\\
&=\prod_{j=1}^m \frac{1}{z^{\lceil D_j\cdot d+ \Lambda_j \D\rceil}} \frac{\Gamma(U_j/z +\langle D_j\cdot d+ \Lambda_j \D\rangle)}{\Gamma(U_j/z+D_j\cdot d+ \Lambda_j \D+1)}.
\end{align*}
The main result of \cite{Jiang2017} states that $I^\p$ lies on the Lagrangian cone for $\p$. Hence $I^\p$ can be used to determine genus $0$ Gromov-Witten theory of $\p$. See \cite{Brown2014} for the case of toric bundles.

Recall the age for $d\in \K$ is defined as $\iota_d=\sum_{i\notin I_d} \langle D_i\cdot d\rangle$, and if $i\in I_d$, we have $D_i\cdot d\in \Z$ so $\langle D_i \cdot d\rangle=0$. So 
\begin{align*}
    \prod_{j=1}^m \frac{1}{z^{\lceil D_j \cdot d+ \Lambda_j \D\rceil}} &= {z^{-\sum_{j=1}^m \lceil D_j\cdot d+ \Lambda_j \D\rceil }} = {z^{-(\sum_{j=1}^m \lceil D_j\cdot d \rceil + \Lambda_j \D )}} \\
    &= {z^{-(\sum_{i=1}^m (D_j\cdot d +\langle D_j\cdot -d\rangle) + (\sum_{i=1}^m\Lambda_j) \D )}}= \frac{1}{z^{\sum_{i=1}^m D_j\cdot d - c_1(P) \D }}\cdot \frac{1}{z^{\iota_{[-d]}}},
\end{align*}
and the $I$-function can be written as
\begin{align*}
I^{\p} &=ze^{\sigma/z} \sum_{\D\in \overline{NE}(B)} J^B_{\D}(z,\tau) Q^{\D} \sum_{d\in\K} y^d \Delta \mathbf{1}_{\left [ -d \right ]}\\
&= ze^{\sigma/z} \sum_{\D\in \overline{NE}(B)} J^B_{\D}(z,\tau) Q^{\D} \sum_{d\in\K}   \frac{y^d}{z^{(D_1+\cdots+D_m) d- c_1(P)\D}}\prod_{j=1}^m \frac{\Gamma(1+U_j/z -\langle- D_j\cdot d\rangle)}{\Gamma(1+U_j/z+D_j\cdot d+ \Lambda_j \D)} \frac{\mathbf{1}_{\left [ -d \right ]}}{z^{\iota_{\left [ -d \right ]}}}.
\end{align*}

Define the $H$-function by defining its degree-$\D$ part to be
\begin{align*}
H^\D &= e^{\sigma/2\pi i} \sum_{d\in \K} y^d \prod_{j=1}^m \frac{1}{\Gamma(1+U_j/2\pi i+D_j\cdot d+ \Lambda_j \D)} \mathbf{1}_{\left[ d\right]}.
\end{align*}
Then we have the following relation between the $H$-function and $I$-function: 
\begin{align}\label{eqn:I_vs_H}
z^{-1} I^\p = \sum_{\D\in \overline{NE}(B)} J^B_{\D}(z,\tau) Q^{\D} {z^{\frac{c_0(\lambda)}{z}}}z^{c_1(P)\D} z^{-\frac{\tx{deg}}{2}} z^{\rho} (\widehat{\Gamma} \cup (2\pi i)^{\frac{\text{deg}_0}{2}} \text{inv}^* H^{\D}(z^{\frac{-\text{deg} (y)}{2}} y)),
\end{align}
where $\widehat{\Gamma}$ is the $T$-equivariant Gamma class given by
$$\widehat{\Gamma}=\bigoplus_{f\in \K/\Ll} \left( \prod_{j=1}^m \Gamma(1+U_j-\langle D_j\cdot f \rangle) \right) \mathbf{1}_f, $$
and $\rho=\theta(D_1+\cdots+D_m)$. The notation $\tx{deg}$ is the degree operator for Chen-Ruan orbifold cohomology and $\tx{deg}_0$ is the degree operator for ordinary cohomology, i.e. the degree without age shifting. The degrees for variables are given by $\tx{deg}(z)=2$ and 
$$\sum_{i=1}^r (\tx{deg}(y_i))p_i = 2\sum_{i=1}^m D_i.$$
By $z^{\frac{-\text{deg} (y)}{2}} y$ we mean
$$(z^{\frac{-\text{deg} (y_1)}{2}} y_1,\dots,z^{\frac{-\text{deg} (y_r)}{2}} y_r).$$

\subsection{Analytic Continuation of $H$-function via Localization}
For a minimal anticone $\delta\in\A$, there is a corresponding $T$-fixed point in each fiber $X_\omega$ over $b\in B$, and together they form a gerbe $B_\delta$ over $B$. Let $i_\delta: B_\delta \ra \p$ be the inclusion map. Similarly, for inertia stack $I\p$, a $T$-fixed gerbe corresponds to a pair $(\delta,f)$ where $\delta\in\A$ and $f\in \K/\Ll$ satisfy $D_i\cdot f\in\Z$ for all $i\in\delta$. Let $i_{(\delta,f)}: B_{(\delta,f)}\ra \p_f$ be the inclusion to $I\p$. We have the following restriction of the $H$-function to a fixed gerbe:
$$i^*_{(\delta,f)} H^\D= \sum_{d\in \K, \left[ d\right]=f} \frac{y^d}{\prod_{j\in \delta} \Gamma(1+D_j\cdot d+ \Lambda_j \D)} \frac{e^{\frac{\sigma(\delta)}{2\pi i}}}{\prod_{j\notin \delta} \Gamma(1+\frac{U_j(\delta)}{2\pi i}+D_j\cdot d+\Lambda_j \D)},$$
where we write $\sigma(\delta)=i^*_{\delta}\sigma$ and $U_j(\delta)=i^*_{\delta}U_j$ for short.

For the rest of this Section, any notation with index $+$ (respectively $-$) means it is for $\p_+$ (respectively $\p_-$). Especially, let $I_\pm$ denote their $I$-functions and $H^\D_\pm$ denote their $H$-functions respectively. $(\delta_+,f_+)$ indexes a $T$-fixed gerbe in $I\p_+$ and $(\delta_-,f_-)$ indexes a $T$-fixed gerbe in $I\p_-$. By $\delta_+|\delta_-$ we mean $\delta_+=\{j_1,\dots,j_{r-1},j_+\}, \delta_-=\{j_1,\dots,j_{r-1},j_-\}$ with $D_{j_1},\dots,D_{j_{r-1}} \in W$, $D_{j_+}\cdot e>0, D_{j_-}\cdot e<0$. By $(\delta_+,f_+)|(\delta_-,f_-)$ we mean $\delta_+|\delta_-$ and there exists $\alpha\in \Q$ such that $f_-=f_+ +\alpha e$.


We will show the following relation between $H^\D_+$ and $H^\D_-$:
\begin{theorem}\label{THM}
Let $(\delta_+,f_+)$ index a T-fixed section on $I\p_+$, then for $\delta_+\in \A_+\cap \A_-$, we have 
$$i^*_{(\delta_+,f_+)} H^\D_+=i^*_{(\delta_+,f_+)} H^\D_-.$$
Otherwise we have 
$$i^*_{(\delta_+,f_+)} H^\D_+=\sum_{(\delta_-,f_-)|(\delta_+,f_+)} C^{\delta_-,f_-}_{\delta_+,f_+} i^*_{(\delta_-,f_-)} H^\D_-,$$
where
\begin{align*}
C^{\delta_-,f_-}_{\delta_+,f_+} &= e^{\frac{\pi i w}{D_{j_-}\cdot e}(\frac{U_{j_-}(\delta_+)}{2\pi i} +D_{j_-} (f_+-f_-)) }\\
&\cdot \frac{\sin \pi(\frac{U_{j_-}(\delta_+)}{2\pi i} +D_{j_-} (f_+-f_-))}{(-D_{j_-}\cdot e)\sin \frac{\pi}{-D_{j_-}\cdot e}(\frac{U_{j_-}(\delta_+)}{2\pi i} +D_{j_-} (f_+-f_-)) } \cdot
\prod_{j:D_j\cdot e<0,j\neq j_-} \frac{\sin\pi(\frac{U_j(\delta_+)}{2\pi i} +D_j\cdot f_+)}{\sin \pi(\frac{U_j(\delta_-)}{2\pi i}+D_j\cdot f_-)},
\end{align*}
with $w=-1-\sum_{j:D_j\cdot e<0} D_j \cdot e=-1 +\sum_{j:D_j\cdot e>0} D_j\cdot e$ and $j_- \in \delta_-$ is the unique element such that $D_{j_-} \cdot e<0$. Note that $C^{\delta_-,f_-}_{\delta_+,f_+}$ does not depend on $\D$.
\end{theorem}

We need the following
\begin{lemma}\label{lem:div_relations}
Let $\delta_-\in \A_-$ and $\delta_+\in \A_+$ be minimal anticones such that $\delta_-|\delta_+$, and $j_-\in\delta_-$ the unique element which is not in $\delta_+$. Then, for any $j\in[m]$, we have
$$U_j(\delta_+)=U_j(\delta_-)+\frac{D_j\cdot e}{D_{j_-}\cdot e} U_{j_-}(\delta_+).$$
\end{lemma}

\begin{proof}
Write $\delta_-=\{j_{1},\dots,j_{r-1},j_-\}$, then $D_1,\dots,D_{r-1},D_{j_-}$ form a basis of $\Ll\otimes \Q$, so for any $D_j$ we have 
\begin{equation}\label{eqn:D_relation}
D_j=c_1 D_1+\cdots+ c_{r-1}D_{r-1} + c_{j_-} D_{j_-}.
\end{equation}
Multiplying (\ref{eqn:D_relation}) with $e$ gives $D_j \cdot e=c_{j_-} D_{j_-}\cdot e$, or $c_{j_-}=\frac{D_j\cdot e}{D_{j_-}\cdot e}$.

Now apply $\theta$ to (\ref{eqn:D_relation}), we have
$$U_j-\lambda_j=c_1 (U_1-\lambda_1)+\cdots+ c_{r-1} (U_{r-1}-\lambda_{r-1})+ c_{j_-} (U_{j_-}-\lambda_{j_-}).$$
Applying $i^*_{\delta_-}$ and $i^*_{\delta_+}$ gives following equations:
$$U_j(\delta_-)-\lambda_j=\sum_{i=1}^{r-1}c_i (-\lambda_i) + c_{j_-} (-\lambda_{j_-}),$$
$$U_j(\delta_+)-\lambda_j=\sum_{i=1}^{r-1}c_i (-\lambda_i) + c_{j_-} (U_{j_-}(\delta_+)-\lambda_{j_-}).$$
Taking their difference, we get $U_j(\delta_+)-U_j(\delta_-)=c_{j_-} U_{j_-}(\delta_+)$, which is what we need.
\end{proof}

\begin{remark}
The above proof shows that this lemma can be generalized to $i^*_{\delta_+}\theta_+(p)=i^*_{\delta_-}\theta_-(p)+\frac{p\cdot e}{D_{j_-}\cdot e} U_{j_-}(\delta_+)$ for any $p\in \Ll\otimes \C$. Especially, we have $\sigma_+(\delta)=\sigma_-(\delta)$ for $\delta\in \A_-\cap\A_+$ and $\sigma_+(\delta_+)=\sigma_-(\delta_-)+\frac{(\tx{log} (y^+))^e}{D_{j_-}\cdot e} U_{j_-}(\delta_+)$ for $\delta_+|\delta_-$.
\end{remark}

\begin{proof}[Proof of Theorem \ref{THM}]
For $\delta\in \A_+\cap \A_-$, the birational map $\varphi:\p_+ \dashrightarrow \p_-$ is an isomorphism in a neighborhood of $B_{\delta}$, and we have $\sigma_+(\delta)=\sigma_-(\delta)$, so $i^*_{(\delta_+,f_+)} H^\D_+=i^*_{(\delta_+,f_+)} H^\D_-$. 

Consider $\delta_+\in \A_+\setminus\A_-$, restricting $H^\D_+$ to $(\delta_+,f_+)$ gives
$$i^*_{(\delta_+,f_+)} H^\D_+= \sum_{d\in \K_+, \left[ d\right]=f_+} \frac{(y^+)^d}{\prod_{j\in \delta_+} \Gamma(1+D_j\cdot d+ \Lambda_j \D)} \frac{e^{\frac{\sigma_+(\delta_+)}{2\pi i}}}{\prod_{j\notin \delta_+} \Gamma(1+\frac{U_j(\delta_+)}{2\pi i}+D_j\cdot d+\Lambda_j \D)}.$$
Define $\delta^{\vee}_+=\{d\in \Ll\otimes \R: D_j\cdot d+ \Lambda_j \D\geq 0, \forall j\in \delta_+\}$. Note that the summation in $i^*_{(\delta_+,f_+)} H^\D_+$ is actually over $d\in \delta^{\vee}_+$. Also for $j\in \delta_+$, $D_j \cdot e=0$ except for one $j_+$ where $D_{j_+}\cdot e> 0$, so for any $d\in \delta^{\vee}_+$ we may write $d=d_+ + ke$ where $d_+\in \delta^{\vee}_+$ but $d_+ -e\notin \delta^{\vee}_+$, $k\in \Z_{\geq0}$. Then 
\begin{align*}
i^*_{(\delta_+,f_+)} H^\D_+= \sum_{d_+\in \delta^{\vee}_+, \left[ d_+\right]=f_+} (y^+)^{d_+} \sum_{k=0}^\infty \frac{e^{\frac{\sigma_+(\delta_+)}{2\pi i}}((y^+)^e)^k}{\prod_{j=1}^m \Gamma(1+\frac{U_j(\delta_+)}{2\pi i}+D_j\cdot d_+ + k D_j\cdot e+\Lambda_j \D)}.
\end{align*}

Using the identity $\Gamma(y)\Gamma(1-y)=\pi/\sin(\pi y)$ to rewrite the second sum over $k$, we have:
\begin{align*}
&\sum_{k=0}^\infty \frac{e^{\frac{\sigma_+(\delta_+)}{2\pi i}}((y^+)^e)^k}{\prod_{j=1}^m \Gamma(1+\frac{U_j(\delta_+)}{2\pi i}+D_j\cdot d_+ + k D_j\cdot e+\Lambda_j \D)} \\
= &\sum_{k=0}^\infty e^{\frac{\sigma_+(\delta_+)}{2\pi i}}((y^+)^e)^k \prod _{j:D_j\cdot e<0} \frac{(-1)^{kD_j\cdot e+\Lambda_j\D} \sin \pi( - \frac{U_j(\delta_+)}{2\pi i}-D_j\cdot d_+)}{\pi} \\
&\times \frac{\prod_{j:D_j\cdot e<0} \Gamma(-(\frac{U_j(\delta_+)}{2\pi i}+D_j\cdot d_+ + k D_j\cdot e+ \Lambda_j \D))}{\prod_{j:D_j\cdot e\geq0} \Gamma(1+\frac{U_j(\delta_+)}{2\pi i}+D_j\cdot d_+ + k D_j \cdot e+ \Lambda_j \D)},
\end{align*}

which is equal to

\begin{align*}
\sum_{k=0}^\infty e^{\frac{\sigma_+(\delta_+)}{2\pi i}}  \text{Res}_{s=k} \Gamma(s)\Gamma(1-s)e^{\pi i s} ((y^+)^e)^s \prod _{j:D_j\cdot e<0} \frac{e^{\pi i(sD_j\cdot e+ \Lambda_j\D)} \sin \pi( - \frac{U_j(\delta_+)}{2\pi i}-D_j\cdot d_+)}{\pi}  \\
\times \frac{\prod_{j:D_j\cdot e<0} \Gamma(-(\frac{U_j(\delta_+)}{2\pi i}+D_j\cdot d_+ + s D_j\cdot e+ \Lambda_j\D))}{\prod_{j:D_j\cdot e\geq0} \Gamma(1+\frac{U_j(\delta_+)}{2\pi i}+D_j\cdot d_+ + s D_j \cdot e+\Lambda_j \D)} ds.
\end{align*}

Consider the following integral
\begin{align*}
e^{\frac{\sigma_+(\delta_+)}{2\pi i}}  \int_{C} \Gamma(s)\Gamma(1-s)
\frac{\prod_{j:D_j\cdot e<0} \Gamma(-(\frac{U_j(\delta_+)}{2\pi i}+D_j\cdot d_+ + s D_j\cdot e+ \Lambda_j\D))}{\prod_{j:D_j\cdot e\geq0} \Gamma(1+\frac{U_j(\delta_+)}{2\pi i}+D_j\cdot d_+ + s D_j \cdot e+\Lambda_j \D)} (e^{-\pi i \omega } (y^+)^e)^s ds.
\end{align*}
It has poles at $s\in\Z$ and 
\begin{equation}\label{eqn:nontrivial_poles}
s=\frac{1}{-D_{j_-}\cdot e} \left(\frac{U_{j_-}(\delta_+)}{2\pi i} +D_{j_-}\cdot d_+ +\Lambda_{j_-} \D -n\right),
\end{equation}
 where $D_{j_-}\cdot e<0$ and $n$ is a non-negative integer. Choose a contour $C$ such that poles at non-negative integers are on the right of $C$ and other poles are on the left of $C$. In view of \cite{Borisov2006} and the assumption $\sum_j D_j\in W$, we can see that this integral is convergent and analytic as a function of $(y^+)^e$ in the domain $\{(y^+)^e: \abs{\tx{arg}((y^+)^e)-w\pi}<\pi\}$. Moreover, let $\cc= \prod_{j:D_j\cdot e\neq 0} (D_j\cdot e )^{D_j\cdot e}$ be the so-called conifold point, then for $\abs{(y^+)^e}<\abs{\cc}$ the integral is equal to the sum of residues on the right of $C$ and for $\abs{(y^+)^e}>\abs{\cc}$ the integral is equal to the negative of the sum of residues on the left of $C$. By our choice of $C$, for $\abs{(y^+)^e}<\abs{\cc}$ the integral, which is equal to the sum of residues on the right of $C$, gives part of the localization of $H^\D_+$, and we define the analytic continuation of $H^\D_+$ using the sum of residues on the left of $C$ for $\abs{(y^+)^e}>\abs{\cc}$. 

To compute the sum of residues on the left, first note that the residues at poles $s=-1-n $ vanish, so we only need to consider poles (\ref{eqn:nontrivial_poles}). By straightforward computations, this is
\begin{equation}\label{eqn:residue}
\begin{split}
&-e^{\frac{\sigma_+(\delta_+)}{2\pi i}} ((y^+)^e)^s  e^{\pi i [s (1+D_{j_-}\cdot e )+ \Lambda_{j_-}\D]} \frac{\sin \pi(\frac{U_{j_-}(\delta_+)}{2\pi i} +D_{j_-}\cdot d_+)}{\sin \pi s} \frac{1}{-D_{j_-}\cdot e} \frac{(-1)^n}{n!}\\
&\times \prod_{j:D_j\cdot e<0,j\neq j_-} \frac{e^{\pi i (sD_j\cdot e+ \Lambda_j \D)} \sin\pi(\frac{U_j(\delta_+)}{2\pi i} +D_j\cdot d_+)}{\sin \pi(\frac{U_j(\delta_+)}{2\pi i}+D_j\cdot d_+ + s D_j\cdot e+\Lambda_j \D)}\\
&\times \prod_{j:j\neq j_-} \frac{1}{\Gamma(1+\frac{U_j(\delta_+)}{2\pi i}+D_j\cdot d_+ + s D_j\cdot e+ \Lambda_j \D)}.
\end{split}
\end{equation}

Since $D_{j_-}\cdot e<0$, $-D_{j_-}\cdot e\in \Z_+$, we may write the non-negative integer $n$ as $n =k(-D_{j_-}e) + \Lambda_{j_-}\D+l$ with $0\leq l+ \Lambda_{j_-}\D<(-D_{j_-}e)$. Let $d_-=d_+ + \frac{D_{j_-}d_+ -l}{-D_{j_-}e} e$, then 
$$-\Lambda_{j_-}\D \leq l= D_{j_-}\cdot d_- <-D_{j_-}e- \Lambda_{j_-}\D, $$
so we have
$$0\leq D_{j_-}\cdot d_- + \Lambda_{j_-}\D < -D_{j_-}\cdot e.$$
Recall our definition of $\delta^{\vee}_+=\{d\in \Ll\otimes \R: D_j\cdot d+ \Lambda_j \D\geq 0, \forall j\in \delta_+\}$, we have $\delta^{\vee}_- =\{d\in \Ll\otimes \R: D_j\cdot d+ \Lambda_j \D\geq 0, \forall j\in \delta_-\}$. So $d_-\in \delta^{\vee}_-$ but $d_- +e \notin\delta^{\vee}_- $. 

By Lemma \ref{lem:div_relations}, we also have 
$$\frac{U_j(\delta_+)}{2\pi i}+D_j\cdot d_+ + s (D_j \cdot e)=\frac{U_j(\delta_-)}{2\pi i}+D_j\cdot d_- - k (D_j\cdot e),$$
and
$$\sigma_+(\delta_+)=\sigma_-(\delta_-)+\frac{(\tx{log} (y^+))^e}{D_{j_-}\cdot e} U_{j_-}(\delta_+).$$

Now we rewrite the residue (\ref{eqn:residue}) to relate it to $i^\star_{(\delta_-,f_-)} H^\D_-$ using these relations:
\begin{equation}\label{eqn:residue2}
\begin{split}
&-e^{\frac{\sigma_-(\delta_-)}{2\pi i}} ((y^+)^{d_- -d_+ -ke})  e^{\frac{\pi i w}{D_{j_-}\cdot e}(\frac{U_{j_-}(\delta_+)}{2\pi i} +D_{j_-} \cdot (d_+-d_-)) } \frac{\sin \pi(\frac{U_{j_-}(\delta_+)}{2\pi i} +D_{j_-}\cdot (d_+-d_-))}{(-D_{j_-}\cdot e)\sin \frac{\pi}{-D_{j_-}\cdot e}(\frac{U_{j_-}(\delta_+)}{2\pi i} +D_{j_-}\cdot (d_+-d_-)) } \\
&\times \prod_{j:D_j e<0,j\neq j_-} \frac{\sin\pi(\frac{U_j(\delta_+)}{2\pi i} +D_j \cdot d_+)}{\sin \pi(\frac{U_j(\delta_-)}{2\pi i}+D_j\cdot d_-)}
\prod_{j=1}^m \frac{1}{\Gamma(1+\frac{U_j(\delta_-)}{2\pi i}+D_j\cdot d_- - k (D_j\cdot e)+\Lambda_j \D)}.
\end{split}
\end{equation}

Let $f_-$ be the equivalent class of $d_-$ in $\K_-/\Ll$, then we have the relation $d_+ - d_-=f_+-f_- +Ne$ for some $N\in\Z$. Take this substitution into (\ref{eqn:residue2}) and the dependence on $N$ cancels, so (\ref{eqn:residue2}) is in fact equal to
$$-e^{\frac{\sigma_-(\delta_-)}{2\pi i}} ((y^+)^{d_- -d_+ -ke}) C^{\delta_-,f_-}_{\delta_+,f_+}\prod_{j=1}^m \frac{1}{\Gamma(1+\frac{U_j(\delta_-)}{2\pi i}+D_j\cdot d_- - k (D_j\cdot e)+\Lambda_j \D)}.$$

Now apply the above result to analytic continuation of $i^*_{(\delta_+,f_+)} H^\D_+$. Note that the summation over residues on the left of $C$ is the same as the summation over $(\delta_-,f_-)$ such that $(\delta_-,f_-)|(\delta_+,f_+)$. We have
\begin{align*}
&i^*_{(\delta_+,f_+)} H^\D_+ = \sum_{d_+\in \delta^{\vee}_+, \left[ d_+\right]=f_+} (y^+)^{d_+} \sum_{k=0}^\infty \frac{e^{\frac{\sigma_+(\delta_+)}{2\pi i}}((y^+)^e)^k}{\prod_{j=1}^m \Gamma(1+\frac{U_j(\delta_+)}{2\pi i}+D_j\cdot d_+ + k (D_j\cdot e)+\Lambda_j \D)}\\
&= \sum_{(\delta_-,f_-)|(\delta_+,f_+)}  \sum_{d_-\in \delta^{\vee}_-, \left[ d_-\right]=f_-} (y^-)^{d_-} \sum_{k=0}^\infty \frac{e^{\frac{\sigma_-(\delta_-)}{2\pi i}}((y^-)^e)^{-k}}{\prod_{j=1}^m \Gamma(1+\frac{U_j(\delta_-)}{2\pi i}+D_j\cdot d_- - k (D_j\cdot e)+\Lambda_j\D)} C^{\delta_-,f_-}_{\delta_+,f_+} \\
&=\sum_{(\delta_-,f_-)|(\delta_+,f_+)} C^{\delta_-,f_-}_{\delta_+,f_+} i^*_{(\delta_-,f_-)} H^\D_-.
\end{align*}
This completes the proof.
\end{proof}

Define the following linear transformation: for $\alpha\in H_{\text{CR},T}^\bullet(\mathcal{P}_+)$, 
\begin{align*}
\U_H(\alpha)=\sum_{(\delta,f):\delta\in \A_+\cap\A_-} (i^*_{(\delta,f)} \alpha) \frac{\mathbf{1}_{\delta,f}}{e_T(N_{\delta,f})}
+ \sum_{(\delta_+,f_+):\delta_+\in \A_+ \setminus \A_-}\sum_{(\delta_-,f_-)|(\delta_+,f_+)} C^{\delta_-,f_-}_{\delta_+,f_+} (i^*_{(\delta_-,f_-)} \alpha) \frac{\mathbf{1}_{\delta_+,f_+}}{e_T(N_{\delta_+,f_+})}.
    \end{align*}
Then it is straightforward to check that $H^\D_+=\U_H H^\D_-$ for all $\D\in \overline{NE}(B) $.

Define a transformation $\U$ to be the one making the following diagram commute:
\begin{equation}\label{diagram:Uh_U}
    \begin{tikzcd}
    H^{\bullet\bullet}_T(I\p_-)\otimes_{\widehat{R}_T} \widehat{S}_T \arrow[r,"\U_H"]\arrow[d,"\tilde{\Psi}_{-}"] & H^{\bullet\bullet}_T(I\p_+)\otimes_{\widehat{R}_T} \widehat{S}_T \arrow[d,"\tilde{\Psi}_{+}"] \\
    \Tilde{\h}(\p_-):=H^{\bullet}_{CR,T}(\p_-)\otimes_{R_T} S_T[\tx{log} z]((z^{-1/k})) \arrow[r,"\U"] & H^{\bullet}_{CR,T}(\p_+)\otimes_{R_T} S_T[\tx{log} z]((z^{-1/k}))=:\Tilde{\h}(\p_+).
    \end{tikzcd}
\end{equation}
Here $\tilde{\Psi}_\pm$ is defined\footnote{$\Gamma(T_B)$ is the Gamma class of the tangent bundle $T_B$.} by the relation (\ref{eqn:I_vs_H}) between $H$-function and $I$-function: 
\begin{equation}\label{eqn:psi_map}
\tilde{\Psi}_\pm(\alpha)=z^{-\frac{\tx{deg}}{2}} z^{\rho_\pm} ({\Gamma(T_B)}\widehat{\Gamma_\pm} \cup (2\pi i)^{\frac{\text{deg}_0}{2}} \text{inv}^* \alpha)|_{Q^{\D}\mapsto Q^\D z^{c_1(\p_\pm)\cdot\D} }. 
\end{equation}
A direct calculation shows $I_+=\U I_-$. This proves part (1) of Theorem \ref{thm:main}. 

Part (2) of Theorem \ref{thm:main} also follows immediately. Recall that $p^+_i=p^-_i,1\leq i\leq r-1$ form a basis for the wall $W$, so part (2) of Theorem \ref{thm:main} is equivalent to $\theta_+(p^+_i)\circ \U = \U\circ \theta_-(p^-_i)$ for $1\leq i\leq r-1$, which follows from definitions of $\theta$ and $\U$.

\subsection{Fourier-Mukai Transformation}
Consider the following diagram: 
\begin{center}
    \begin{tikzcd}
     & \tilde{\p}\arrow[dl,"f_+"']\arrow[dr,"f_-"] &  \\
    \p_+ \arrow[rr,dotted,"\varphi"] & & \p_-,
    \end{tikzcd}
\end{center}
where $\p_+, \p_-$ are defined as above and $\tilde{\p}$ is the toric stack bundle with fiber $\tilde{X}$ which is the toric blow-up of $X_+, X_-$, and $f_+,f_-$ are induced by the blow-up maps of the fibers. The Fourier-Mukai functor $R(f_+)_*\circ Lf_-^*: D^b_T(\p_-)\to D^b_T(\p_+)$ is an equivalence of category, see \cite{CCT}. We consider the isomorphism $\fm:K^0_T(\p_-)\ra K^0_T(\p_+)$ induced by this Fourier-Mukai equivalence.

Recall that the fibers $X_\pm$ of the toric stack bundle $\p_\pm$ are defined from characters $D_1,\dots,D_m \in \Ll^\vee$, stability conditions $\omega_\pm \in \Ll^\vee\otimes\R$ and we have the short exact sequence:
$$0\ra \Ll \overset{(D_1,\dots,D_m)}{\xrightarrow{\hspace{1cm}}} \Z^m \overset{\beta}{\longrightarrow} \N \ra 0. $$
The toric blow-up $\Tilde{X}$ can be defined from $\tilde{D}_1,\dots,\Tilde{D}_{m+1}$ and a stability condition $\tilde{\omega}$. Consider the following short exact sequence:
$$0\ra \Ll\oplus\Z \overset{(\tilde{D}_1,\dots,\tilde{D}_{m+1})}{\xrightarrow{\hspace{2cm}}} \Z^{m+1} \overset{\tilde{\beta}}{\longrightarrow} \N \ra 0, $$
where $\tilde{D}_i$ is defined by
\begin{align*}
    \tilde{D_i}=\begin{cases}
 D_i\oplus 0 &\text{ if } 1\leq j\leq m \ \tx{and}\ D_i\cdot e\leq 0, \\ 
  D_i\oplus (-D_i\cdot e) &\text{ if } 1\leq j\leq m \ \tx{and}\ D_i\cdot e> 0, \\ 
  0\oplus 1 &\text{ if } j=m+1. 
\end{cases}
\end{align*}
Let $C_\pm$ be the chamber containing $\omega_\pm$ and $W$ the wall between the two chambers. The stability condition $\tilde{\omega}\in (\Ll\oplus\Z)^\vee\otimes\R$ is defined as $(\omega_0,-\epsilon)$ where $\omega_0$ is in the relative interior of $W\cap \overline{C_+}=W\cap \overline{C_-}$ and $\epsilon$ is a small positive real number. By \cite[Lemma 6.16]{Coates2018}, the toric stack $\tilde{X}$ defined from the above data is the common toric blow-up of $X_\pm$. 

A calculation of the $K$-theory of toric stack bundles can be found in \cite{JT}. We first give a brief description of $K^0_T(\p)$ for a toric stack bundle $\p$. Given an irreducible $K$-representation $p\in \tx{Hom}(K,\C^\times)=\Ll^\vee$, we define a line bundle $L(p)$ over $\p$ as follow:
$$L(p)=(P\times_{(\C^\times)^m} U_\omega)\times \C / (b,z,v)\sim (b,gz,g(p)v),$$
where $g\in K, b\in B, z\in U_\omega, v\in \C$. This line bundle defines a class in $K^0_T(\p)$. Especially, let $p=D_i$, we define the following class in $K^0_T(\p_{\pm})$:
$$R^{\pm}_i=L_{\pm} (D_i) \otimes  e^{\lambda_i},$$ 
where $\pm$ stands for line bundles over $\p_{\pm}$ respectively, $e^{\lambda_i}\in\C[T]$ stands for the irreducible $T$-representation given by the $i$-th projection $T\ra \C^\times$. Note that these line bundles correspond to the toric divisors, and their inverses are denoted by $S^{\pm}_i$. 

Similarly, we may define a line bundle $L(p,n)$ over $\tilde{\p}$ for a character $(p,n)\in \tx{Hom}(K\times \C^\times ,\C^\times)=\Ll^\vee\oplus \Z$, and classes $\tilde{R}_i=L (\tilde{D}_i) \otimes  e^{\lambda_i} $ for $1\leq i\leq m$ and $\tilde{R}_{m+1}=L(\tilde{D}_{m+1})=L(0,1)$ in $K^0_T(\tilde{\p})$. Note that we have $f^*_-(L_-(p))=L(p,0)$ and $f^*_+(L_+(p))=L(p,-p\cdot e)$. If we let $k_i=\tx{max}(D_i\cdot e,0)$ and $l_i=\tx{max}(-D_i\cdot e,0)$, then we have $f^*_-R^-_i=\tilde{R}_i\tilde{R}_{m+1}^{k_i}$ and $f^*_+R^+_i=\tilde{R}_i\tilde{R}_{m+1}^{l_i}$. 

Recall that a $T$-fixed locus, which is a gerbe over $B$, is indexed by a minimal anticone $\delta\in\A$. Let $i_\delta: B_\delta \ra \p$ be the inclusion of this fixed locus. Let $G_\delta$ be the common isotropy group of $ B_\delta$ and $\varrho\in \tx{Hom}(G_\delta,\C^\times)$ an irreducible representation, then choose a lift $\hat{\varrho}\in \tx{Hom}(K,\C^\times)=\Ll$. A localized basis for $K^0_T(\p)$ over $K^0(B)\otimes \C[T]$ is given by
$$\left\{e_{\delta,\varrho}:=(i_\delta)_* \varrho=L(\hat{\varrho})\prod_{i\notin \delta} (1-S_i)\right\}.$$ 
We use $e_{\delta_{\pm},\varrho}$ to denote elements of the basis for $\p_{\pm}$ respectively, and we will describe the Fourier-Mukai transformation on the basis.\footnote{$\mathbb{FM}$ is identity on classes pulled back from $B$.}

Consider a minimal anticone $\tilde{\delta}=\{j_1,\dots,j_{r-1},j_-,j_+\}\in \A_{\tilde{\omega}}$ for $\tilde{\p}$ with all $j_1, \dots ,j_{r-1}, j_-, j_+\in[m]$, $D_{j_1}\cdot e=\cdots=D_{j_{r-1}}\cdot e=0$, $D_{j_-}\cdot e<0$ and $D_{j_+}\cdot e>0$. Let ${\delta_-}=\{j_1,\dots,j_{r-1},j_-\}\in \A_{-}$, ${\delta_+}=\{j_1,\dots,j_{r-1},j_+\}\in \A_{+}$ be the corresponding anticones for $\p_-$ and $\p_+$. Then we have the following diagram:
\begin{center}
    \begin{tikzcd}
      B_{\delta_-}\arrow[r,"i_{\delta_-}"] & \p_- \\ B_{\tilde{\delta}}\arrow[u,"f_{-,\tilde{\delta}}"']\arrow[d,"f_{+,\tilde{\delta}}"]\arrow[r,"i_{\tilde{\delta}}"] & \tilde{\p} \arrow[u,"f_-"']\arrow[d,"f_+"]\\
      B_{\delta_+}\arrow[r,"i_{\delta_+}"] & \p_+.
    \end{tikzcd}
\end{center}

\begin{lemma}\label{lem:pull-push_bundles}
For a Laurent polynomial $q$, we have
$$(f_{+,\tilde{\delta}})_*(i_{\tilde{\delta}})^* L(p,n) q(\tilde{R_1},\dots,\tilde{R}_{m+1}) = (i_{\delta_+})^* \frac{1}{l} \sum_{t\in \mathcal{T}} L_+(p) t^{p\cdot e + n} q(t^{-l_1}R^+_1,\dots,t^{-l_m}R^+_m,t), $$
where $l=-D_{j_-}\cdot e $ and $\mathcal{T}=\{ \zeta\cdot (R^+_{j_-})^{1/l}: \zeta\in \boldsymbol{\mu}_l \}$.
\end{lemma}

\begin{proof}
    The common isotropy groups for $B_{\tilde{\delta}}$ and $B_{\delta_+}$ are given by
    \begin{align*}
        G_{\tilde{\delta}} &= \{ (g,\lambda)\in K\times \C^\times: D_{j_-}(g)=1, D_j(g) \lambda^{-D_j\cdot e }=1 \ \tx{for} \ j\in \delta_+ \},\\
        G_{\delta_+} &= \{ h\in K: D_j(h)=1 \ \tx{for all} \ j\in\delta_+  \}.
    \end{align*}
    The map $f_{+,\tilde{\delta}}$ is induced by a homomorphism from $G_{\tilde{\delta}}$ to $G_{\delta_+}$, $\phi_+: (g,\lambda)\mapsto h=g\lambda^{-e}$. The kernel of this map is $\{(\lambda^e,\lambda): \lambda\in \boldsymbol{\mu}_l \}$, here $\boldsymbol{\mu}_l$ is the cyclic group of order $l$. This gives the following exact sequence
    $$1\ra \boldsymbol{\mu}_l \ra G_{\tilde{\delta}} \ra G_{\delta_+} \ra 1. $$
    This actually shows that $f_{+,\tilde{\delta}}$ makes $B_{\tilde{\delta}}$ a $\boldsymbol{\mu}_l$-gerbe over $B_{\delta_+}$. Now given a character $(p,n)\in \tx{Hom}(K\times \C^\times,\C^\times)$, recall that the line bundle $L(p,n)$ is defined by
    $$L(p,n) =\left( (\tilde{P}\times_{(\C^\times)^{m+1}} U_{\tilde{\omega}})\times \C \right) / (b,z,v)\sim (b,(g,\lambda)\cdot z,p(g)\lambda^n v),$$
    where $\tilde{P}$ is a principal $(\C^\times)^{m+1}$-bundle, $U_{\tilde{\omega}}$ is defined by the GIT data $\tilde{\omega}$, $\tilde{D}_1,\dots,\tilde{D}_{m+1}$ and $\A_{\tilde{\omega}}$ as introduced in Section \ref{CTS}, and for any $(g,\lambda)\in K\times \C^\times$, $z=(z_1,\dots,z_{m+1})\in U_{\tilde{\omega}}$, the action is defined by $\tilde{D}_1,\dots,\tilde{D}_{m+1}$, especially we have $$(g,\lambda)\cdot z = (\tilde{D}_1(g,\lambda) z_1,\dots, \tilde{D}_m(g,\lambda) z_m, \lambda z_{m+1}).$$
    
    First we consider $(i_{\tilde{\delta}})^* L(p,n)$. On each fiber, the line bundle is pulled back to the fixed point corresponding to $\tilde{\delta}=\tilde{\delta}=\{j_1,\dots,j_{r-1},j_-,j_+\}$:
    $$\{ (z_1,\dots,z_{m+1})\in U_{\tilde{\omega}}: z_i=0\ \tx{if}\ i\notin \tilde{\delta} \}/(K\times \C^\times).$$
    In other words, the possibly non-zero coordinates are $(z_{j_1},\dots,z_{j_{r-1}},z_{j_-},z_{j_+})$, and the action of $(g,\lambda)$ on these coordinates is given by
    $$(D_{j_1}(g)z_{j_1},\dots,D_{j_{r-1}}(g)z_{j_{r-1}},D_{j_-}(g)z_{j_-},D_{j_+}(g)\lambda^{-D_{j_+}\cdot e} z_{j_+}).$$
    Now take $(g,\lambda)=(\lambda^e,\lambda),\lambda\in \boldsymbol{\mu}_l $ and consider the equivalence relation in $L(p,n)$, we have $(b,z,v)\sim (b,z,\lambda^{p\cdot e +n} v)$ when restricting to the fixed point, so it vanishes if $l$ does not divide $p\cdot e +n$. If $l$ divides $p\cdot e +n$, first we have
    \begin{align*}
        (f_{+,\tilde{\delta}})^* (i_{\delta_+})^* L_+(p) (R^+_{j_-})^{(p\cdot e+n)/l} &= (i_{\tilde{\delta}})^* (f_+)^* L_+(p) (R^+_{j_-})^{(p\cdot e+n)/l} \\
        &= (i_{\tilde{\delta}})^* L(p,-p\cdot e) (\tilde{R}_{j_-})^{(p\cdot e+n)/l} (\tilde{R}_{m+1})^{p\cdot e +n}\\
        &= (i_{\tilde{\delta}})^* L(p,-p\cdot e) (\tilde{R}_{m+1})^{p\cdot e +n}\\
        &= (i_{\tilde{\delta}})^* L(p,n).
    \end{align*}
    Applying $(f_{+,\tilde{\delta}})_*$ to the above, we get
    $$(f_{+,\tilde{\delta}})_*(i_{\tilde{\delta}})^* L(p,n) = (i_{\delta_+})^* L_+(p) (R^+_{j_-})^{(p\cdot e+n)/l}.$$
We can write the two cases together as 
    $$(f_{+,\tilde{\delta}})_*(i_{\tilde{\delta}})^* L(p,n) = (i_{\delta_+})^* \frac{1}{l} \sum_{t\in\T} L_+(p) t^{p\cdot e+n},$$
    where $\mathcal{T}=\{ \zeta\cdot (R^+_{j_-})^{1/l}: \zeta\in \boldsymbol{\mu}_l \}$: this is because 
    \begin{align*}
      \frac{1}{l} \sum_{t\in\T} t^n =  \left\{\begin{matrix}
       (R^+_{j_-})^{n/l} \ &\tx{if}\ l\ \tx{divides}\ n,\\
       0         &\tx{otherwise}.
\end{matrix}\right.
    \end{align*}
    To prove the lemma, it suffices to consider the case of monomials $q(\tilde{R}_1,\dots,\tilde{R}_{m+1})= \prod_{i=1}^{m+1} \tilde{R}_i^{r_i}$:
    \begin{align*}
        &(f_{+,\tilde{\delta}})_*(i_{\tilde{\delta}})^* L(p,n) q(\tilde{R_1},\dots,\tilde{R}_{m+1}) \\
        &= (f_{+,\tilde{\delta}})_*(i_{\tilde{\delta}})^* L(p,n) \prod_{i=1}^{m+1} \tilde{R}_i^{r_i}\\
        &= (f_{+,\tilde{\delta}})_*(i_{\tilde{\delta}})^* L(p+\sum_{i=1}^m r_i D_i,n+r_{m+1}-\sum_{i=1}^m r_ik_i)\otimes e^{\sum_{i=1}^m r_i \lambda_i} \\
        &= (i_{\delta_+})^* \frac{1}{l} \sum_{t\in\T} L_+(p+\sum_{i=1}^m r_i D_i) t^{p\cdot e+ n+r_{m+1}-\sum_{i=1}^m r_ik_i} \otimes e^{\sum_{i=1}^m r_i \lambda_i}\\
        &=(i_{\delta_+})^* \frac{1}{l} \sum_{t\in\T} L_+(p) t^{p\cdot e+ n} \prod_{i=1}^m (t^{-k_i})^{r_i} (R^+_i)^{r_i} t^{r_{m+1}} \\
        &= (i_{\delta_+})^* \frac{1}{l} \sum_{t\in \mathcal{T}} L_+(p) t^{p\cdot e + n} q(t^{-l_1}R^+_1,\dots,t^{-l_m}R^+_m,t).
    \end{align*}
\end{proof}

\begin{theorem}[]
Let $\delta\in \A_-\cap \A_+$ be a common minimal anticone, then we have
$$\fm(e_{\delta,\varrho})=e_{\delta,\varrho},$$
where $\delta$ is regarded as a minimal anticone in $\A_-$ on the left and in $\A_+$ on the right. Otherwise, for $\delta_-\in \A_-\setminus \A_+$ we have
$$\fm(e_{\delta_-,\varrho})=\frac{1}{l} \sum_{t\in\T} \frac{1-S^+_{j_-}}{1-t^{-1}} L_+(\hat{\varrho}) t^{\hat{\varrho}\cdot e} \prod_{i\notin \delta_-} (1- t^{-k_i} S^+_i).$$
\end{theorem}

\begin{proof}
For $\delta\in \A_+\cap \A_-$, the birational map $\varphi:\p_+\dashrightarrow\p_-$ is an isomorphism in a neighborhood of $B_{\delta}$, so $\fm(e_{\delta,\varrho})=e_{\delta,\varrho}$. Now for $\delta_-\in \A_-\setminus \A_+$, let $\delta_+\in \A_+\setminus \A_-$ be such that $\delta_+|\delta_-$. Consider a minimal anticone $\tilde{\delta}=\{j_1,\dots,j_{r-1},j_-,j_+\}\in \A_{~}$ for $\tilde{\p}$ with all $\{j_1,\dots,j_{r-1},j_-,j_+\}\in[m]$, $D_{j_1}\cdot e=\cdots=D_{j_{r-1}}\cdot e=0$, $D_{j_-}\cdot e<0$ and $D_{j_+}\cdot e>0$. Let ${\delta_-}=\{j_1,\dots,j_{r-1},j_-\}\in \A_{-}$, ${\delta_+}=\{j_1,\dots,j_{r-1},j_+\}\in \A_{+}$ be the corresponding anticones for $\p_-$ and $\p_+$. Recall the following diagram:
\begin{center}
    \begin{tikzcd}
      B_{\delta_-}\arrow[r,"i_{\delta_-}"] & \p_- \\ B_{\tilde{\delta}}\arrow[u,"f_{-,\tilde{\delta}}"']\arrow[d,"f_{+,\tilde{\delta}}"]\arrow[r,"i_{\tilde{\delta}}"] & \tilde{\p} \arrow[u,"f_-"']\arrow[d,"f_+"]\\
      B_{\delta_+}\arrow[r,"i_{\delta_+}"] & \p_+.
    \end{tikzcd}
\end{center}
Take $e_{\delta_{-},\varrho}\in K^0_T(\p_-)$, then we have $(f_-)^* e_{\delta_{-},\varrho}=(f_-)^* L_-(\hat{\varrho})\prod_{i\notin \delta_-} (1-S^-_i)=L(\hat{\varrho},0)\prod_{i\notin \delta_-} (1-\tilde{S}_i\tilde{S}_{m+1}^{k_i})$. By localization formula, we have 
\begin{align*}
(f_+)_*(f_-)^* e_{\delta_{-},\varrho} &= (f_+)_* L(\hat{\varrho},0)\prod_{i\notin \delta_-} (1-\tilde{S}_i\tilde{S}_{m+1}^{k_i}) \\
 &= (f_+)_* \sum_{\tilde{\delta}:\delta_-\subset \tilde{\delta}} (i_{\tilde{\delta}})_* (i_{\tilde{\delta}})^* 
 \frac{L(\hat{\varrho},0)\prod_{i\notin \delta_-} (1-\tilde{S}_i\tilde{S}_{m+1}^{k_i})}{(1-\tilde{S}_{m+1})(\prod_{i\notin \tilde{\delta}} 1-\tilde{S}_i)}\\
 &= \sum_{\tilde{\delta}:\delta_-\subset \tilde{\delta}} (i_{\delta_+})_* (f_{+,\tilde{\delta}})_* (i_{\tilde{\delta}})^*  \frac{L(\hat{\varrho},0)\prod_{i\notin \delta_-} (1-\tilde{S}_i\tilde{S}_{m+1}^{k_i})}{(1-\tilde{S}_{m+1})(\prod_{i\notin \tilde{\delta}} 1-\tilde{S}_i)}.
\end{align*}
Applying Lemma \ref{lem:pull-push_bundles}, we have
 \begin{align*}
 (f_+)_*(f_-)^* e_{\delta_{-},\varrho} &= \sum_{\delta_+: \delta_+ | \delta_-} (i_{\delta_+})_* (i_{\delta_+})^* (\frac{1}{l} \sum_{t\in\T} L_+(\hat{\varrho}) t^{\hat{\varrho}\cdot e} \frac{\prod_{i\notin \delta_-} (1- t^{l_i-k_i} S^+_i)}{(1-t^{-1}) \prod_{i\notin \tilde{\delta}} (1- t^{l_i} S^+_i )})\\
 &= \sum_{\delta_+: \delta_+ | \delta_-} (i_{\delta_+})_* (i_{\delta_+})^* \frac{\frac{1}{l} \sum_{t\in\T}  L_+(\hat{\varrho}) t^{\hat{\varrho}\cdot e} \frac{1-S^+_{j_-}}{1-t^{-1}} \prod_{i\notin \delta_-} (1- t^{-k_i} S^+_i)} {\prod_{i\notin \delta_+} (1- S^+_i)}.
\end{align*}
We need to check that restriction of
$$\frac{1}{l} \sum_{t\in\T}  L_+(\hat{\varrho}) t^{\hat{\varrho}\cdot e} \frac{1-S^+_{j_-}}{1-t^{-1}} \prod_{i\notin \delta_-} (1- t^{-k_i} S^+_i)$$
is only non-zero for fixed points corresponding to minimal anticones $\delta_+$ next to $\delta_-$. 

First recall that $(i_\delta)^* S^+_i=1$ for $i\in \delta$, so for a minimal anticone $\delta$, if there is $i\in \delta$ but $i\notin \delta_-$ such that $D_i\cdot e\leq 0$, then the restriction of $\prod_{i\notin \delta_-} (1- t^{-k_i} S^+_i)$ vanishes, so the restriction of the above summation is $0$. Another case is that $j_-\in \delta$ and there is some $j_+\in \delta$ such that $D_{j_+}\cdot e>0$. Then for each $t\in\T$, there is a factor
$$\frac{1-S^+_{j_-}}{1-t^{-1}}(1- t^{-D_{j_+}\cdot e} S^+_{j_+})= (1-S^+_{j_-})\frac{1-(t^{-1})^{D_{j_+}\cdot e}}{1-t^{-1}}.$$
Since $D_{j_+}\cdot e$ is a positive integer, this factor is a polynomial in $t^{-1}$ and its restriction to $\delta$ vanishes. So the only possibility for non-zero restriction occurs for $\delta_+$ with $\delta_+|\delta_-$. Applying localization formula one more time, we obtain
 $$\fm(e_{\delta_-,\varrho})= \frac{1}{l} \sum_{t\in\T} \frac{1-S^+_{j_-}}{1-t^{-1}} L_+(\hat{\varrho}) t^{\hat{\varrho}\cdot e} \prod_{i\notin \delta_-} (1- t^{-k_i} S^+_i).$$
\end{proof}

\subsection{Comparison}
We claim that the Fourier-Mukai transformation is compatible with the analytic continuation above, i.e. the following diagram commutes
\begin{equation}\label{diagram:FM_Uh}
\begin{tikzcd}
K^0_T(\p_-) \arrow[r,"\fm"]\arrow[d,"\tilde{\tx{ch}}"] & K^0_T(\p_+) \arrow[d,"\tilde{\tx{ch}}"] \\
H^{\bullet\bullet}_T(I\p_-) \arrow[r,"\U_H"] & H^{\bullet\bullet}_T(I\p_+). 
\end{tikzcd}
\end{equation}
Here $\tilde{\tx{ch}}$ is the orbifold Chern character.

First, we check commutativity on basis $K^0_T(B)$. Note that $\tilde{\tx{ch}}$ maps $K^0_T(B)$ to $H^{\bullet}(B)\otimes R_T$ and everything is fixed under $\fm$ and $\U_H$ in this case.

Now for a basis element $e_{\delta_-,\varrho}\in K^0_T(\p_-)$, if $\delta_-\in \A_-\cap \A_+$, we also have
$$\tilde{\tx{ch}}(\fm(e_{\delta_{-},\varrho}))=\tilde{\tx{ch}}(e_{\delta_{-},\varrho})=\U_H \tilde{\tx{ch}}(e_{\delta_{-},\varrho}).$$
For $\delta_-\notin \A_-\cap \A_+$, we calculate them explicitly. First, note that we have
\begin{align*}
\tilde{\tx{ch}}(L_{\pm}(\hat{\varrho})) &= \bigoplus_{f\in \K_{\pm}/\Ll} e^{2\pi i \hat{\varrho}\cdot f } e^{\theta_{\pm}(\hat{\varrho})} \mathbf{1}_f.
\end{align*}
Especially we have
\begin{align*}
\tilde{\tx{ch}}(S^{\pm}_j) &= \bigoplus_{f\in \K_{\pm}/\Ll} e^{-2\pi i D_j\cdot f } e^{-U_j} \mathbf{1}_f,
\end{align*}
and
\begin{align*}
\tilde{\tx{ch}}(t) &= \bigoplus_{f\in \K_{\pm}/\Ll} \zeta e^{2\pi i D_{j_-}\cdot f /l} e^{U_{j_-}/l} \mathbf{1}_f,
\end{align*}
for $t=\zeta (R^+_{j_-})^{1/l} \in \T$. 
We may rewrite $\fm(e_{\delta_-,\varrho})$ as
$$\fm(e_{\delta_-,\varrho})= \sum_{t\in\T} \left(\frac{1-S^+_{j_-}}{l(1-t^{-1})} \prod_{i\notin \delta_-, D_i\cdot e<0} \frac{1- S^+_i}{1- t^{-D_i \cdot e} S^+_i}\right) \left(L_+(\hat{\varrho}) t^{\hat{\varrho}\cdot e} \prod_{i\notin \delta_-} (1- t^{-D_i \cdot e} S^+_i)\right).$$
It suffices to check the equality locally. Note that $\tilde{\tx{ch}}(\fm(e_{\delta_{-},\varrho}))$ is only supported on fixed loci corresponding to $(\delta_+,f_+)$ where $\delta_+| \delta_-$. So we need to check
$$i^\star _{(\delta_+,f_+)} \tilde{\tx{ch}}(\fm(e_{\delta_{-},\varrho})) = \sum_{(\delta'_-,f_-)|(\delta_+,f_+)} C^{\delta'_-,f_-}_{\delta_+,f_+} (i^\star _{(\delta'_-,f_-)} \tilde{\tx{ch}}(e_{\delta_{-},\varrho})).$$
First, using Lemma \ref{lem:div_relations} we have 
\begin{align*}
i^\star_{(\delta_+,f_+)} \tilde{\tx{ch}}(L_+(\hat{\varrho}) t^{\hat{\varrho}\cdot e}) &= e^{2\pi i \hat{\varrho}\cdot f_+ } e^{i_{\delta_+}^\star \theta_{+}(\hat{\varrho})} e^{-2\pi i D_{j_-}\cdot f_- \frac{\hat{\varrho}\cdot e}{l}} e^{2\pi i D_{j_-}\cdot f_+ \frac{\hat{\varrho}\cdot e}{l}} e^{U_{j_-}(\delta_+)\frac{\hat{\varrho}\cdot e}{l}}\\
&= e^{2\pi i \hat{\varrho}\cdot f_+} e^{2\pi i \frac{\hat{\varrho}\cdot e}{l} D_{j_-}\cdot(f_+-f_-)} e^{i^\star_{\delta_+} \theta_+(\hat{\varrho})+ \frac{\hat{\varrho}\cdot e}{l} U_{j_-}(\delta_+)}\\
&= e^{2\pi i (\hat{\varrho}\cdot f_+ - \alpha \hat{\varrho}\cdot e)} e^{i^\star_{\delta_-} \theta_-(\hat{\varrho})}\\
&= e^{2\pi i \hat{\varrho}\cdot f_-}e^{i^\star_{\delta_-} \theta_-(\hat{\varrho})}\\
&= i^\star_{(\delta_-,f_-)} \tilde{\tx{ch}}(L_-(\hat{\varrho})),
\end{align*}
where we set $f_+=f_- + \alpha e$ and $t=e^{-2\pi i D_{j_-}\cdot f_- /l} (R^+_{j_-})^{1/l}$. Also
\begin{align*}
i^\star_{(\delta_+,f_+)} \tilde{\tx{ch}}(S^+_j t^{-D_j \cdot e}) &= e^{-2\pi i D_j\cdot f_+} e^{-U_j (\delta_+)} e^{-2\pi i D_{j_-}\cdot f_- \frac{-D_j\cdot e}{l}} e^{2\pi i D_{j_-}\cdot f_+ \frac{-D_j\cdot e}{l}} e^{U_{j_-}(\delta_+)\frac{-D_j\cdot e}{l}}\\
&= e^{-2\pi i \left(D_j\cdot f_+ - \frac{D_j\cdot e}{D_{j_-}\cdot e} D_{j_-}(f_+-f_-)  \right)} e^{-\left(U_j(\delta_+)- \frac{D_j\cdot e}{D_{j_-}\cdot e} U_{j_-}(\delta_+)\right)}\\
&= e^{-2\pi i D_j\cdot f_-} e^{-U_j(\delta_-)}\\
&= i^\star_{(\delta_-,f_-)} \tilde{\tx{ch}}(S_j^-).
\end{align*}
Together we have
\begin{align*}
i^\star_{(\delta_+,f_+)} \tilde{\tx{ch}}\left(L_+(\hat{\varrho}) t^{\hat{\varrho}\cdot e} \prod_{i\notin \delta_-} (1-S^+_i t^{-D_i \cdot e}) \right) = i^\star_{(\delta_-,f_-)} \tilde{\tx{ch}}(e_{\delta_-,\varrho}).
\end{align*}

On the other hand, we have
\begin{align*}
& i^\star_{(\delta_+,f_+)} \tilde{\tx{ch}}\left(\frac{1-S^+_{j_-}}{l(1-t^{-1})} \prod_{i\notin \delta_-, D_i\cdot e<0} \frac{1- S^+_i}{1- t^{-D_i \cdot e} S^+_i}\right)\\
&= \frac{1}{l} \frac{1-e^{-U_{j_-}(\delta_+)-2\pi i D_{j_-}\cdot f_+}}{1- e^{-\frac{1}{l}(U_{j_-}(\delta_+)+2\pi i D_{j_-}(f_+ -f_-))}} \prod_{i\notin \delta_-, D_i\cdot e<0} \frac{1-e^{-U_i(\delta_+)-2\pi i D_i\cdot f_+}}{1-e^{-U_i(\delta_-)-2\pi i D_i\cdot f_-} }\\
&= \frac{1}{l} \frac{\sin\left( \frac{-U_{j_-}(\delta_+)}{2i}-\pi  D_{j_-}\cdot f_+ \right)}{\sin\left( -\frac{1}{l}(\frac{U_{j_-}(\delta_+)}{2i}+\pi  D_{j_-}(f_+ -f_-))\right)} \prod_{i\notin \delta_-, D_i\cdot e<0} \frac{\sin\left( \frac{-U_i(\delta_+)}{2i}-\pi D_i\cdot f_+\right)} { \sin\left( \frac{-U_i(\delta_-)}{2i}-\pi D_i\cdot f_- \right) }\\
&\times e^{-U_{j_-}(\delta_+)/2- \pi i D_{j_-}\cdot f_+ + \frac{1}{2l}(U_{j_-}(\delta_+)+ 2\pi i D_{j_-}(f_+ -f_-)) + \sum_{i\notin \delta_-, D_i\cdot e<0} (-U_i(\delta_+)/2-\pi i D_i\cdot f_+ + U_i(\delta_-)/2+\pi i D_i\cdot f_-) }\\
&= \frac{1}{l} e^{\frac{\pi i w}{D_{j_-}\cdot e} (\frac{U_{j_-}(\delta_+)}{2\pi i}+ D_{j_-}(f_+-f_-)) } 
\frac{\sin \pi(\frac{U_{j_-}(\delta_+)}{2\pi i} +D_{j_-} (f_+-f_-))}{\sin \frac{\pi}{-D_{j_-}e}(\frac{U_{j_-}(\delta_+)}{2\pi i} +D_{j_-} (f_+-f_-)) } 
\prod_{i\notin \delta_-, D_i\cdot e<0} \frac{\sin\pi(\frac{U_i(\delta_+)}{2\pi i} +D_i f_+)}{\sin \pi(\frac{U_i(\delta_-)}{2\pi i}+D_i f_-)}\\
&= C^{\delta_-,f_-}_{\delta_+,f_+}. 
\end{align*}
So for any basis element $e_{\delta_{-},\varrho}\in \A_-\setminus \A_+$ and any $\delta_+| \delta_-$, we have
\begin{align*}
 i^\star_{(\delta_+,f_+)} \tilde{\tx{ch}}(\fm(e_{\delta_{-},\varrho})) &= \sum_{t\in \T} i^\star_{(\delta_+,f_+)} \tilde{\tx{ch}} \left(\frac{1-S^+_{j_-}}{l(1-t^{-1})} \prod_{i\notin \delta_-, D_i\cdot e<0} \frac{1- S^+_i}{1- t^{-D_i \cdot e} S^+_i}\right)\\
 &\times i^\star_{(\delta_+,f_+)} \tilde{\tx{ch}}\left(L_+(\hat{\varrho}) t^{\hat{\varrho}\cdot e} \prod_{i\notin \delta_-} (1- t^{-D_i \cdot e} S^+_i)\right)\\
 &= \sum_{(\delta_-,f_-)|(\delta_+,f_+)} C^{\delta_-,f_-}_{\delta_+,f_+} i^\star_{(\delta_-,f_-)} \tilde{\tx{ch}}(e_{\delta_-,\varrho}) \\
 &= i^\star_{(\delta_+,f_+)} \U_H \tilde{\tx{ch}}(e_{\delta_-,\varrho}),
\end{align*}
since as $f_-$ runs over $\K_-/ \Ll$, $\zeta= e^{-2\pi i D_{j_-}\cdot f_- /l}$ runs over $\boldsymbol{\mu}_l$. This completes the proof of commutativity of (\ref{diagram:FM_Uh}). 

Combining diagrams (\ref{diagram:FM_Uh}) and (\ref{diagram:Uh_U}), we arrive at part (3) of Theorem \ref{thm:main}.

The claim that $\U$ is symplectic follows from the commutativity of (\ref{diagram:FM_compatible}), the fact that $\mathbb{FM}$ preserves Mukai pairing, and Toen's Grothendieck-Riemann-Roch formula.


\bibliographystyle{amsplain}

\end{document}